\documentclass{article}

\pdfoutput=1 
\usepackage[english]{babel}
\usepackage{csquotes}
\usepackage{biblatex}
\usepackage{indentfirst}
\usepackage{authblk}
\usepackage{amsmath}
\usepackage{caption}
\usepackage{amsthm}
\usepackage{booktabs}
\usepackage{color}
\usepackage{tabu}
\usepackage{fancyhdr}
\usepackage{authblk}
\usepackage{amssymb}
\usepackage[letterpaper,top=2cm,bottom=2cm,left=3cm,right=3cm,marginparwidth=1.75cm]{geometry}
\usepackage{mathtools}
\usepackage{graphicx}
\usepackage{fancyhdr}
\usepackage{array}
\usepackage[colorlinks=true, allcolors=blue]{hyperref}
\usepackage{graphicx}
\usepackage{float} 
\usepackage{subfigure}

\title{Koopman Spectral Linearization vs.  Carleman Linearization: \quad A Computational Comparison Study}
\author{Dongwei Shi, Xiu Yang}
\affil{Department of Industrial and System Engineering, Lehigh University}

\begin{document}
\maketitle
\begin{abstract}
Nonlinearity presents a significant challenge in problems involving dynamical systems, prompting the exploration of various linearization techniques, including the well-known Carleman Linearization. In this paper, we introduce the Koopman Spectral Linearization method tailored for nonlinear autonomous dynamical systems. This innovative linearization approach harnesses the Chebyshev differentiation matrix and the Koopman Operator to yield a lifted linear system. It holds the promise of serving as an alternative approach that can be employed in scenarios where Carleman linearization is traditionally 
applied. Numerical experiments demonstrate the effectiveness of this linearization approach for several commonly used nonlinear dynamical systems.
\end{abstract}
\par {\bf\emph{Keywords: Dynamical systems; Koopman Operator; Carleman Linearization.}}

\section{Introduction}\label{sec1}
Nonlinear dynamical systems find extensive applications in both science and engineering. However, analyzing and computing solutions for nonlinear systems can pose a formidable challenge. A natural strategy for mitigating this complexity involves approximating a nonlinear system through a technique known as linearization. The most conventional method of linearization entails employing a first-order Taylor expansion of the system dynamics. While this technique has proven effective and is widely used in Model Predictive Control (MPC) and similar domains, it remains confined to a very narrow vicinity around the expansion point. To comprehensively address these limitations, a series of more sophisticated techniques and theories have emerged. Approaches such as Carleman linearization and the Koopman operator are rooted in the concept of lifting nonlinear systems to infinite-dimensional linear systems. These global linearization techniques and theories were developed by scholars like Torsten Carleman \cite{ref10}, Bernard Koopman \cite{ref11}, and John Von Neumann \cite{ref12} in the 1930s. In practice, by employing finite section approximation, we can attain a finite yet higher-dimensional linear system. This empowers us to capture the characteristics of the original nonlinear system over a longer time period compared with the first-order Taylor expansion approach.

\par Within the control system community and quantum computing community, the application of Carleman linearization concepts has yielded a range of notable accomplishments across various studies. For instance, in \cite{ref1}, the author emphasizes the applicability of Carleman Linearization in designing nonlinear controllers and state estimators. The work in \cite{ref2} harnessed Carleman approximation to establish a correlation between the lifted system and the domain of attraction of the original nonlinear system. Noteworthy recent endeavors, including \cite{ref3}, skillfully harnessed Carleman linearization to implement efficient model predictive control mechanisms for nonlinear systems. The application of Carleman linearization concepts was extended to \cite{ref4}, where they were adeptly utilized for tasks such as state estimation and the formulation of feedback control laws. Moreover, \cite{ref5,ref6} exploited the inherent structure of the elevated system to introduce a practical methodology for quantitatively addressing and solving the Hamilton-Jacobi-Bellman equation through a meticulously crafted iterative approach.
Furthermore, in the realm of quantum computing, Carleman linearization has also found prominence. Leveraging its ability to transform nonlinear dynamical systems into linear forms, researchers have employed Carleman linearization to design efficient quantum algorithms, with provable expoential speed up with respect to dimension compared with best possible classic algorithms \cite{ref7,ref8,ref9}.

\par
As highlighted in reference \cite{ref13}, the Carleman matrix is, in fact, a transposed finite section matrix approximation of the Koopman Operator Generator in the polynomial basis. Despite this connection, modern research on the Koopman Operator and Carleman Linearization follows distinct paths. Investigations into the Koopman Operator predominantly adopt a data-driven approach, as they use  spatiotemporal data without requiring the explicitly equation to find an approximation of the principal Koopman modes and correspondingly Koopman eigenpairs. Notably numerical schemes like Dynamic Mode Decomposition (DMD) \cite{ref18,ref19,ref20}  and its diverse variations, including Extended Dynamic Mode Decomposition (EDMD) \cite{ref21}, make use of data collected from time snapshots of a specified dynamical system. These approaches employ Singular Value Decomposition (SVD) and other matrix decomposition techniques to reveal the temporal and spatial spectral attributes of the system. Expanding beyond DMD, EDMD and their variants, methods for obtaining tractable representations have been further connected with deep neural networks. \cite{ref22,ref23} Two specific neural network architectures have been explored in this context. One is based on the autoencoder principle, employing a low-dimensional latent space to enhance the interpretability of outcomes. The other architecture involves a higher-dimensional latent space, which often yields improved results when dealing with dynamical systems featuring continuous spectra. A comprehensive review of the theoretical framework and algorithms related to Data-Driven Koopman Operators is provided by \cite{ref15}.

\par
Diverging from the previously discussed data-driven methods, our approach constitutes a progressive evolution of the Linearization Technique. This distinction arises from our specific need for an analytical equation. More precisely, our method integrates the Chebyshev Differentiation Matrix in the spectral method \cite{ref24,ref25}, yielding a clear matrix approximation of the Koopman Operator's generator. Subsequently, this matrix is employed to generate an explicit higher-dimensional linear equation that closely approximates the original nonlinear system. In this context, the identical approximation matrix was employed across various scalar observable configurations, with only the initial value being altered. This innovative approach was motivated by the numerical ordinary differential equation (ODE) algorithm discussed in \cite{ref16,ref17}. This study centers around autonomous systems and potentially serves as a viable alternative approach applicable in situations where Carleman linearization has conventionally been utilized.

\par
The paper is structured as follows: \autoref{sec2} presents the background knowledge, followed by the discussion of the Koopman Spectral Linearization Method in \autoref{sec3}. In \autoref{sec4}, we provide the numerical results, and the subsequent sections delve into the discussion and conclusions in \autoref{sec5}.

\section{Background}\label{sec2}
\subsection{Carleman Linearization}\label{sec2.1}
Carleman Linearization provides a systematic methodology for deriving the corresponding lifted dynamics from an initial value problem defined as follows \cite{ref27}: 
\begin{equation}
\begin{aligned}
    \frac{d\boldsymbol{x}}{dt}
    &=B_1\boldsymbol{x}+B_2\boldsymbol{x}^{\otimes2}+...+B_k\boldsymbol{x}^{\otimes k}\quad\text{with}\quad \boldsymbol{x}_0=\boldsymbol{x}(t_0)
\label{eq1}
\end{aligned}
\end{equation}
Here the matrices $B_i\in \mathcal{R}^{d\times d^i}$ are time-independent, k is the degree of the polynomial ODE and $\boldsymbol{x}^{\otimes i}$  
denotes the j-fold tensor product for each positive integer j. We define an infinite-dimensional vector $\boldsymbol{y}=(\boldsymbol{y_1},\boldsymbol{y_2},\boldsymbol{y_3},....)^T$ with $ \boldsymbol{y_i}=\boldsymbol{x}^{\otimes i}$, and formulate matrices $A_{i+j-1}^i$
\begin{equation}
\begin{aligned}
A_{i+j-1}^i=\sum_{v=1}^{i}\otimes_iB_j
\end{aligned}
\label{eq2}
\end{equation}
Where $\otimes_iB=I\otimes I\otimes I \otimes ...\otimes B\otimes ...\otimes I$ with $B$ at $(d-i+1)$-th position.
Further, formula \eqref{eq2} lead to the following equation
\begin{equation}
\begin{aligned}
\frac{d\boldsymbol{y_i}}{dt}=\sum_{j=1}^{k}A_{i+j-1}^i\boldsymbol{y_{i+j-1}}
\end{aligned}
\label{eq3}
\end{equation}
Formally, this can be represented as an infinite-dimensional ordinary differential equation:
\begin{equation}
\begin{aligned}
\frac{d\boldsymbol{y}}{dt}&=\boldsymbol{\mathcal{A}}\boldsymbol{y} \qquad \text{with}\quad
\boldsymbol{y_i}(t_0)=\boldsymbol{x}(t_0)^{\otimes i}
\end{aligned}
\label{eq4}
\end{equation}
where $\boldsymbol{\mathcal{A}}$ is an infinite-dimensional matrix defined as:
\begin{equation}
\begin{aligned}
\boldsymbol{\mathcal{A}}=
\begin{bmatrix}
    &A_1^1 &A_2^1&A_3^1 &...&A_k^1&0&0&...\\
    &0&A_2^2&A_3^2&...&A_k^2&A_{k+1}^2&0&...\\
    &0&0&A_3^3 &...&A_k^3&A_{k+1}^3&A_{k+2}^3&...\\
    &\vdots &\vdots&\vdots&&\vdots&\vdots&\vdots&\\
\end{bmatrix}
\end{aligned}
\label{eq5}
\end{equation}
In practice, the system is often truncated at a certain order 
N. Here is an illustrative example \cite{ref15} of Carleman Linearization applied to the following system:
\begin{equation}
\begin{aligned}
\frac{dx}{dt}&=x^2
\end{aligned}
\label{eq6}
\end{equation}
Carleman Linearization procedure can derive a linear ode like 
\begin{equation}
\begin{aligned}
\frac{d}{dt}
\begin{bmatrix}
    x\\x^{2}\\x^{3}\\ \vdots
\end{bmatrix}
=
\begin{bmatrix}
0&1&0&0&0&...\\
0&0&2&0&0&...\\
0&0&0&3&0&...\\
\vdots&\vdots&\vdots&\vdots&\vdots
\end{bmatrix}
\begin{bmatrix}
    x\\x^{2}\\x^{3}\\ \vdots
\end{bmatrix}
\end{aligned}
\label{eq7}
\end{equation}
\subsection{Koopman Operator}\label{sec2.2}
Firstly consider d-dimension (i.e. $\boldsymbol{x}\in X \subset \mathbb{R}^d$) dynamical system described by a first order autonomous order ordinary differential equations with $t\in [t_0,T]$.

\begin{equation}
\frac{d\boldsymbol{x}(t)}{dt}=\boldsymbol{f}(\boldsymbol{x}(t))  
\label{eq8}
\end{equation}
where $\boldsymbol{x}=(x_1,x_2,...,x_d)^T$ belongs to a space X and this dynamics $\boldsymbol{f}$ is also an d dimensional vector valued function with $(f_1,f_2,...,f_d)^T$ serve as each coordinate. For our analysis, we introduce observables or measurement functions, represented by the function $g:X\rightarrow\mathbb{R}$, note this is a function on $\mathcal{G}(X)$ (i.e. $\mathcal{L}^2$ Hilbert space). And
the flow map $F_t$ represents the evolution of the system dynamics as a mapping on X.

\begin{equation}
F_t:\boldsymbol{x}(t_0)\rightarrow \boldsymbol{x}(t+t_0)
\label{eq9}
\end{equation}
 \\
The Koopman operator family, denoted as $\{\mathcal{K}_t\}$, is defined as:

\begin{equation}
\mathcal{K}_t:g(\boldsymbol{x}(t_0))\rightarrow g(F_t(\boldsymbol{x}(t_0)))=g(\boldsymbol{x}(t_0+t))
\label{eq10}
\end{equation}
Alternatively, we can express $\mathcal{K}_t$ as a composition of functions:

\begin{equation}
\mathcal{K}_t:g\rightarrow g\circ F_t
\label{eq11}
\end{equation}
It's crucial to note that the Koopman operator is linear, owing to the linearity of the function space $\mathcal{G}(X)$. Further we can define the generator of the Koopman operator
$\mathcal{K}$

\begin{equation}
\mathcal{K}g
:=\lim_{t\rightarrow0}\frac{\mathcal{K}_tg-g}{t}
=\frac{\partial g}{\partial t}
\label{eq12}
\end{equation}
The generator can be explicitly formulated as follows:

\begin{equation}
\mathcal{K}g
=\frac{\partial g}{\partial t}
=\nabla g(\boldsymbol{x})\frac{d\boldsymbol{x}}{dt}
=\nabla g(\boldsymbol{x})\cdot \boldsymbol{f}(\boldsymbol{x})
\label{eq13}
\end{equation}
The transformation of representing it in the form of an operator can be further extended as follows:
\begin{equation}
\mathcal{K}
=\nabla \cdot \boldsymbol{f}(\boldsymbol{x})
=\sum_{i=1}^d f_i(x)\frac{\partial}{\partial x_i}
\label{eq14}
\end{equation}
As a trade-off, we convert a nonlinear mapping on the variable $x$ into a linear operator on the variable $g$, which, in turn, leads to an infinite-dimensional space. 

In Koopman Spectral Theory, eigenvalues and eigenfunctions (or eigenpairs for brevity) are employed for practical numerical computations\cite{ref14,ref19}. Suppose $(\lambda,\phi(\boldsymbol{x}))$ is an eigenpair for the generator $\mathcal{K}$. By definition

\begin{equation}
\begin{aligned}
\mathcal{K}\phi(\boldsymbol{x}(t+t_0))
=\frac{ d\phi(\boldsymbol{x}(t+t_0))}{dt}
=\lambda \phi(\boldsymbol{x}(t+t_0))
\end{aligned}
\label{eq15}
\end{equation}
Based on ode in equation (15), we can derive the following outcomes:
\begin{equation}
    \phi(\boldsymbol{x}(t+t_0))=\phi(\boldsymbol{x}(t_0))e^{\lambda t}
\label{eq16}
\end{equation}
Actually, $\phi(\boldsymbol{x})$ is also an eigen-function for $\mathcal{K}_t$,this can be observed from

\begin{equation}
    \mathcal{K}_t\phi(\boldsymbol{x}(t+t_0))=\phi(\boldsymbol{x}(t_0+2t))=e^{\lambda t}\phi(\boldsymbol{x}(t_0+t))
\label{eq17}
\end{equation}
Therefore if we denote $e^{\lambda t}$ as $\mu$ than $(\mu,\phi)=(e^{\lambda t},\phi)$ is an eigenpair for $\mathcal{K}_t$. 

If $g \in \mathcal{G}(X) \subset Span\{\phi_k\}$, then even nonlinear observable can be represented as a linear combination of eigenfunctions. This can be formulated as  $g=\sum_jc_j\phi_j$, where $c_j\in\mathbb{R}$ is the corresponding coefficient of g with respect to $\phi_j$. Further as pointed in \cite{ref14}, the evolution of the observable can also be represented as a linear combination of eigenpairs which is derived from:
\begin{equation}
    \mathcal{K}_tg(\boldsymbol{x}(t+t_0))
    =\mathcal{K}_t(\sum_{j=1}^{\infty}c_j\phi_j(\boldsymbol{x}(t_0+t)))
    =\sum_{j=1}^{\infty}c_je^{\lambda_j t}\phi_j(\boldsymbol{x}(t_0+t))
    =\sum_{j=1}^{\infty}c_j\mu_j\phi_j(\boldsymbol{x}(t_0+t))
    \label{eq18}
\end{equation}
Hence, we can infer
\begin{equation}
    g(\boldsymbol{x}(t+t_0))
    =\sum_{j=1}^{\infty}c_je^{\lambda_j t}\phi_j(\boldsymbol{x}(t_0))
    \label{eq19}
\end{equation}
Multi-dimensional observables are based on same manner. Just consider vector-valued $c_j$  with same dimension of the observable. $c_j$ is also known as j-th Koopman mode \cite{ref15}.

\section{Koopman Spectral Linearization Method}\label{sec3}
\subsection{Construction of the lifted matrix}\label{3.1}
By adopting construction of the matrix approximation from \cite{ref16} and \cite{ref17}, the generator of Koopman Operator can be approximated as following:\\
Currently, when $d$=1, for some eigen-pairs $(\phi,\lambda)$
based on equation \eqref{eq7}
\begin{equation}
\mathcal{K}\phi= f\cdot\frac{\partial}{\partial x}\cdot\phi
\label{eq20}
\end{equation}
To derive a finite dimensional approximation, we starting from polynomial interpolation of the eigenfunction $\phi(x)$ on spatial discretized Guass-Labatto points $\{\xi_i\}_{i=1}^N$. Which give us
\begin{equation}
\phi(x)\approx \phi^N(x)=\sum_{i=1}^{N}\phi^N(\xi_i)L_i(x) 
\label{eq21}
\end{equation}
with $L_j$ are Lagrange Polynomials such that $L_j(\xi_i)=\delta_{ij}$ where $\delta_{ij}$ is the delta function.
 Therefore we obtained
\begin{equation}
K
\begin{bmatrix}
\phi^N(\xi_0)\\\\
\phi^N(\xi_1)\\\\
\vdots\\\\
\phi^N(\xi_N)
\end{bmatrix}
=
	\begin{bmatrix}
	 f(\xi_0) &   & &  \\
	  & f(\xi_1) &  & \\
	  &   & \ddots&\\
      &   &    &f(\xi_N)
	 \end{bmatrix}
  D
\begin{bmatrix}
\phi^N(\xi_0)\\\\
\phi^N(\xi_1)\\\\
\vdots\\\\
\phi^N(\xi_N)
\end{bmatrix}
\label{eq22}
\end{equation}
Based on equation \eqref{eq22}, the finite representation of $K$ is 
\begin{equation}
    K=diag\{f(\xi_0),f(\xi_1),...,f(\xi_N)\}D
\label{eq23}
\end{equation}
When $d=2$, Let  $\{\xi_i\}_{i=1}^N$,and$\{\eta_j\}_{j=1}^N$ be the Guass-Lobatto Points of $x_1$ and $x_2$. And denote $\Theta=\{(\xi_i,\eta_j)\}_{i,j=1}^N$ Now for every bivariate eigenfunction $\phi(x_1,x_2)$ we have polynomial interpolation $\phi^N$
\begin{equation}
    \phi(x_1,x_2)\approx \phi^N(x_1,x_2)=\sum_{i=1}^N\sum_{j=1}^N\phi^N(\xi_i,\eta_j)L_i(x_1)L_j(x_2)
\label{eq24}
\end{equation}
Hence all function values on the collocation points can be formed as matrix denoted as $\phi^N(\boldsymbol{\Theta})$
\begin{equation}
    \phi^N(\boldsymbol{\Theta})=
    \begin{bmatrix}
        \phi^N(\xi_1,\eta_1)&\phi^N(\xi_1,\eta_2)&...&\phi^N(\xi_1,\eta_N)\\\phi^N(\xi_2,\eta_1)&\phi^N(\xi_2,\eta_2)&...&\phi^N(\xi_2,\eta_N)\\
        \vdots&\vdots&\ddots&\vdots\\
        \phi^N(\xi_N,\eta_1)&\phi^N(\xi_N,\eta_2)&...&\phi^N(\xi_N,\eta_N)
    \end{bmatrix}
\label{eq25}
\end{equation}
Let $D_1,D_2$ be the differentiation matrices for $x_1$ and $x_2$ respectively, and $f_1(\boldsymbol{\Theta})$,$f_2(\boldsymbol{\Theta})$ be the matrices of $f_1$ $f_2$ evaluated at $(\xi_i,\eta_j)$. And we denote $\mathcal{K}\phi^N(\boldsymbol{\Theta})_{ij}=\mathcal{K}\phi^N(\xi_i,\eta_j)$ Then the formula below can be derived based on \eqref{eq7}
\begin{equation}
\mathcal{K}\phi^N(\boldsymbol{\Theta})=f_1(\boldsymbol{\Theta})\odot(D_1\phi^N(\boldsymbol{\Theta}))+f_2(\boldsymbol{\Theta})\odot(\phi^N(\boldsymbol{\Theta}) D_2^T)
\label{eq26}
\end{equation}
We vectorize all the matrix operations as vector operations, this process can be formulated as:
\begin{equation}
 \begin{aligned}
 \mathcal{K}vec(\phi^N(\boldsymbol{\Theta}))&=vec(f_1(\boldsymbol{\Theta}))\odot((I\otimes D_1)vec(\phi^N(\boldsymbol{\Theta})))+vec(f_2(\boldsymbol{\Theta}))\odot((D_2\otimes I)vec(\phi^N(\boldsymbol{\Theta})))\\
  &=[diag(vec(f_1(\boldsymbol{\Theta})))(I\otimes D_1)+ diag(vec(f_2(\boldsymbol{\Theta}))(D_2\otimes I))] vec(\phi^N(\boldsymbol{\Theta})\\
  &=Kvec(\phi^N(\boldsymbol{\Theta})
 \end{aligned}
 \label{eq27}
\end{equation}
Therefore, for $d=2$ cases, the construction of matrix $K$ is given by 
\begin{equation}
    K=diag(vec(f_1(\boldsymbol{\Theta})))(I\otimes D_1)+ diag(vec(f_2(\boldsymbol{\Theta}))(D_2\otimes I))
\label{eq28}
\end{equation}
In general for a d-dimensional case, consider Guass-Labotta points for each coordinate as \\$\{\xi_{i_1}^1\}_{i_1=1}^N$,$\{\xi_{i_2}^2\}_{i_2=1}^N$,...,$\{\xi_{i_d}^d\}_{i_d=1}^N$ And use $\Theta$ denotes the collection of all the collocation points i.e. $\Theta=\{(\xi_{i_1}^1,\xi_{i_2}^2,...,\xi_{i_d}^d):i_1,i_2,...,i_d,\text{traverse from 1 to N}\}$. And the eigenfunction $\phi$ is approximated by
\begin{equation}
\phi(x_1,...,x_d)\approx\phi^N(x_1,...,x_d)=\sum^{N}_{i_1,...,i_d=1}\phi^N(\xi_{i_1}^1,\xi_{i_2}^2,...,\xi_{i_d}^d)L_{i_1}(x_1)...L_{i_d}(x_d)
\label{eq29}
\end{equation}
Corresponingly, $\phi^N(\boldsymbol{\Theta})$ is a tensor of eigenfunction values on collocation points. and $D_1,D_2,...D_d$ means differentiation matrices. With n-mode multiplication in tensor algebra, we can write
\begin{equation}
\mathcal{K}\phi^N(\boldsymbol{\Theta})
=\sum_{i=1}^d f_i(\Theta)\odot(\phi^N(\boldsymbol{\Theta)}\times_iD_i)
\label{eq30}
\end{equation}
Pursuing the analogous vectorization approach in equation \eqref{eq27}, we can reformulate equation \eqref{eq30} as
\begin{equation}
\begin{aligned}
\mathcal{K}vec(\phi^N(\boldsymbol{\Theta}))
&=\{\sum_{i=1}^d vec(f_i(\boldsymbol{\Theta}))\odot(\otimes_iD_i)\}vec(\phi^N(\boldsymbol{\Theta}))\\
&=\{\sum_{i=1}^d diag(vec(f_i(\boldsymbol{\Theta})))(\otimes_iD_i)\}vec(\phi^N(\boldsymbol{\Theta}))
\end{aligned}
\label{eq31}
\end{equation}
Consequently,
\begin{equation}
    K=\sum_{i=1}^d diag(vec(f_i(\boldsymbol{\Theta})))(\otimes_iD_i)
\label{eq32}
\end{equation}
Where $\otimes_iD_i=I\otimes I\otimes I \otimes ...\otimes D_i\otimes ...\otimes I$ with $D_i$ at $(d-i+1)$-position, to be more intuitive, here is an example when $d=3$. 
\begin{equation}
\begin{aligned}
K
=&diag\{vec(F_1)\}I\otimes I\otimes D_1\\
+&diag\{vec(F_2)\}I\otimes D_2\otimes I\\
+&diag\{vec(F_3)\}D_3\otimes I\otimes I
\end{aligned}
\label{eq33}
\end{equation}
\subsection{Solution of the linear system}\label{sec3.2}
Now, let's examine the solution of the linear system to provide further insight into why it serves as an approximation of the original nonlinear one. Referring back to equation \eqref{eq19}, we have. \begin{equation}
    g(\boldsymbol{x}(t))\approx g^N(\boldsymbol{x}(t))=\sum_{j=1}^N\hat{c}_j\phi^N_j(\boldsymbol{x}(t_0))e^{\hat{\lambda}_j}
    \label{eq34}
\end{equation} When $d=1$, for a scalar observable g. Consider Guass-Lobatto points on a region with radius r around $x_0$, where $\Theta=\{\xi_i\}_{i=1}^N$ and $\xi_1<\xi_2<...<\xi_N$. Here we pick odd number N such that $\xi_{(N+1)/2}=x_0=x(t_0)$. Thus all Guass-Lobatto Points are on the interval $[x_0-r,x_0+r]$. Suppose $(\hat{\lambda}_j,v_j)$ is an eigen-pair of matrix K. Vector $v_j$ are used to approximate eigenfunction $\phi^N_j$ evaluated at the collocation points. we have $(v_j)_i=\phi^N_j(\xi_i)\approx \phi_j(\xi_i)$. Hence
\begin{equation}
g^N(x(t))=\sum_{j=1}^N\hat{c}_j\phi^N_j(x(t_0))e^{\hat{\lambda}_jt}=\sum_{j=1}^N\hat{c}_j\phi^N_j(\xi_{(N+1)/2})e^{\hat{\lambda}_jt}=
\sum_{j=1}^N\hat{c}_j(v_j)_{\frac{N+1}{2}}e^{\hat{\lambda}_jt}
\label{eq35}
\end{equation}
Now consider the N-dimensional linear system below
\begin{equation}
\frac{dy}{dt}=Ky\qquad y_0=(g(x_0-r),...,g(x_0+r))^T   
\label{eq36}
\end{equation}
We can write the analytical solution
\begin{equation}
    y(t)=e^{Kt}y_0
\label{eq37}
\end{equation}

 Suppose the eigen-decompostion of $K=V\Lambda V^{-1}$, with each columns of V are eigenvectors of K, denoted as $v_j$, and $\Lambda$ contain all the eigenvalues. Further by setting $t=0$ in equation \eqref{eq35}, we have
 \begin{equation}
     g_N(x_0)=\sum_{j=1}^N\hat{c}_j\phi^N_j(x_0)
     \label{eq38}
 \end{equation}
 This is also true for different  initial value thus
 \begin{equation}
     g_N(\xi_i)=\sum_{j=1}^N\hat{c}_j\phi^N_j(\xi_i)=\sum_{j=1}^N\hat{c}_j(v_j)_i,\quad i=1,2,...,N 
     \label{eq39}
 \end{equation}Thus we have $V^{-1}y_0=c$ since $y(0)=(g^N(\xi_1),g^N(\xi_2),...,g^N(\xi_N)^T$, where c is called the Koopman mode.
\begin{equation}
\begin{aligned}
    y(t)
    &=e^{Kt}y_0\\
    &=e^{V\Lambda tV^{-1}}y_0\\
    &=Ve^{\Lambda t}V^{-1}y_0\\
    &=\begin{bmatrix}
        v_1,v_2,...,v_N
    \end{bmatrix}
	\begin{bmatrix}
	 e^{\lambda_1t} &   & &  \\
	  & e^{\lambda_2t} &  & \\
	  &   & \ddots&\\
	  &   &  &e^{\lambda_Nt}
	 \end{bmatrix}
    \begin{bmatrix}
        c_1\\c_2\\\vdots\\c_N
    \end{bmatrix}\\
    &=\sum_{j=1}^Nc_je^{\lambda_jt}v_j
\end{aligned}
\label{eq40}
\end{equation}
And $g^N(x(t))=\sum_{j=1}^Nc_je^{\lambda_jt}(v_j)_{\frac{N}{2}}$. 

Notice $c_j$ are Koopman modes, $e^{\lambda_jt}$ are eigenfunctions.  
In general case, for d-dimensional nonlinear system, suppose $\boldsymbol{\Theta}$ represents d-dimensional collocation points, and K is constructed as equation \eqref{eq32}. Now consider the linear system given below
\begin{equation}
    \frac{d\boldsymbol{y}(t)}{dt}=K\boldsymbol{y} \qquad \boldsymbol{y}(0)=vec(g(\boldsymbol{\Theta}))
    \label{eq41}
\end{equation}
Again by vectorization, $\phi_j(\boldsymbol{x_0})$ can be approximated by the "middle term" of tensor $\phi^N(\boldsymbol{\Theta})$
which leads to 
\begin{equation}
    y(t)=\sum_{j=1}^{N^d}c_je^{\lambda_jt}(v_j)
    \label{eq42}
\end{equation}
And to compute the observable, only the middle term is needed for this long vector.
\begin{equation}
g_N(x(t))=\sum_{j=1}^{N^d}c_je^{\lambda_jt}(v_j)_{\frac{N^d+1}{2}}
\label{eq43}
\end{equation}
which is exactly the solution of \eqref{eq41}.
\par
Unlike Carleman Linearization, which captures the information of all variables within a single linear ordinary differential equation, our linearization approach is primarily tailored for scalar observables. When the need arises to compute multi-dimensional observables or distinct scalar observables, we can still employ the previously constructed lifted matrix $K$. However, this requires transforming the initial conditions accordingly. For instance, when computing an observable $g(\boldsymbol{x})=\boldsymbol{x}=(g_1(\boldsymbol{x}),g_2(\boldsymbol{x}),g_3(\boldsymbol{x}))$,we would need to utilize three different initial conditions: $g_1(\boldsymbol{\Theta}),g_2(\boldsymbol{\Theta}),g_3(\boldsymbol{\Theta})$
This adaptation allows us to extend the applicability of our approach to scenarios involving diverse observable functions since $K$ can be reused.
\par As for the solution of Carleman Linearization, the matrix $\boldsymbol{\mathcal{A}}$ shall be truncated at certain level $N$ which yields $\boldsymbol{\mathcal{A}_N}$ matrix as following:
\begin{equation}
    \boldsymbol{\mathcal{A}_N}=
\begin{bmatrix}
    &A_1^1 &A_2^1&A_3^1 &...&A_k^1&0&0&...\\
    &0&A_2^2&A_3^2&...&A_k^2&A_{k+1}^2&0&...\\
    &0&0&A_3^3 &...&A_k^3&A_{k+1}^3&A_{k+2}^3&...\\
    &\vdots &\vdots&\vdots&&\vdots&\vdots&\vdots&\\
    &0&0&0&....&0&0&0&A_N^N\\
\end{bmatrix}
\end{equation}
Then the infinite linear system can be approximated by 
\begin{equation}
    \frac{dy}{dt}=\boldsymbol{\mathcal{A}_N}y
\end{equation} and the approximated solution of original nonlinear system is exactly the first $d$ elements of the solution for equation (45).
\newpage
\section{Numerical Results}\label{sec4}
In this section, we present the performance of Koopman Spectral Linearization on five nonlinear dynamic models in \autoref{sec4.1}. In each example, we investigate the influence of truncation order $N$ and radius $r$ on accuracy. The reference solution is generated by RK4 if there is no closed-form solution available. To be more precise, in two or three-dimensional cases, we test the influence of the radius in each direction separately. Next, in \autoref{sec4.2}, we compare the accuracy and computational efficiency involved with Carleman Linearization. These comparisons include errors against truncation order, matrix size, and computational time.
\subsection{Linearize Dynamics with Koopman Spectral Linearization}\label{sec4.1}
\subsubsection{Quadratic Model}\label{sec4.1.1}
The quadratic model is possible simplest nonlinear ode just designed for demonstration purpose. The governing ODE is given by
$$
\frac{dx}{dt}=x^2
$$
To be more clarified, we set $x(0)=0.08$ and T=10 in this example. This quadratic model has a closed form solution $x(t)=\frac{1}{(1/x_0)-t}$. The target time period is $[0,T]$. The three numerical experiments use the following parameter settings:
\begin{itemize}
  \item Test of Truncation Order $N$: $r=0.03$
  \item Test of radius $r$: $N=11$
\end{itemize}

\begin{figure}[H]
\centering  
\subfigure[Test of Truncation Order]{
\label{Fig.1a}
\includegraphics[width=0.45\textwidth]{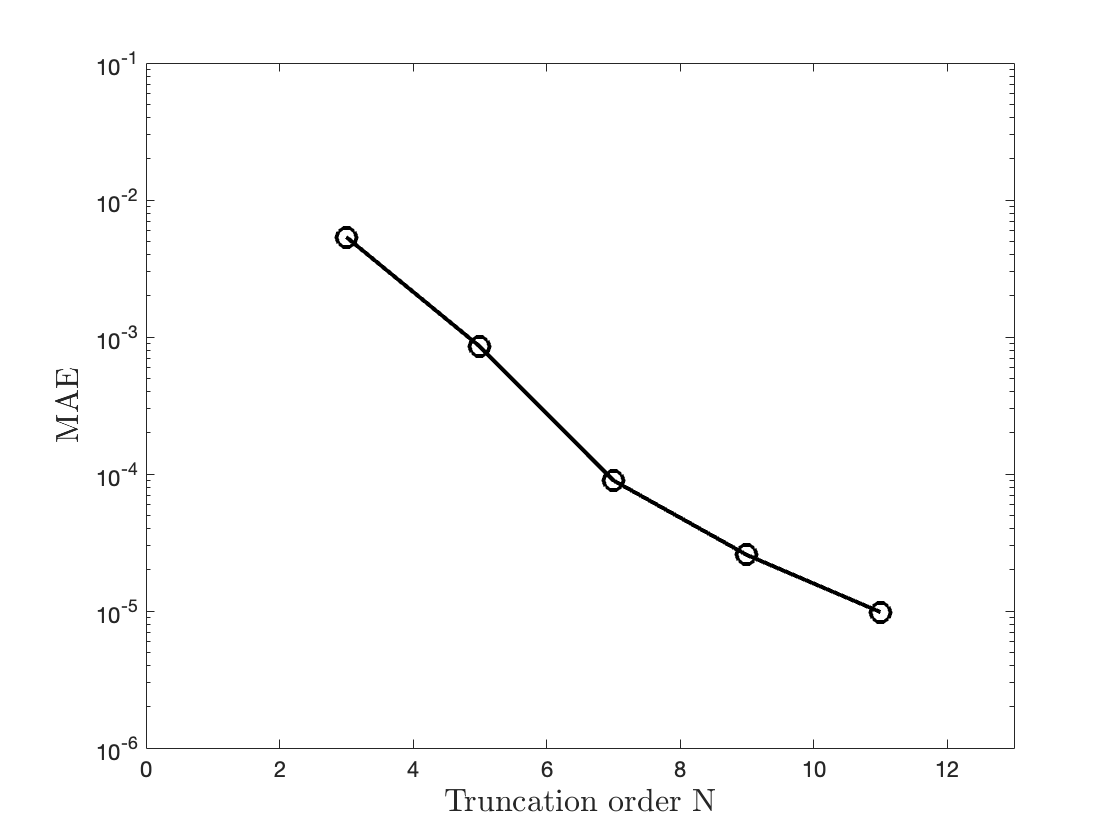}}
\subfigure[Test of Radius]{
\label{Fig.1b}
\includegraphics[width=0.45\textwidth]{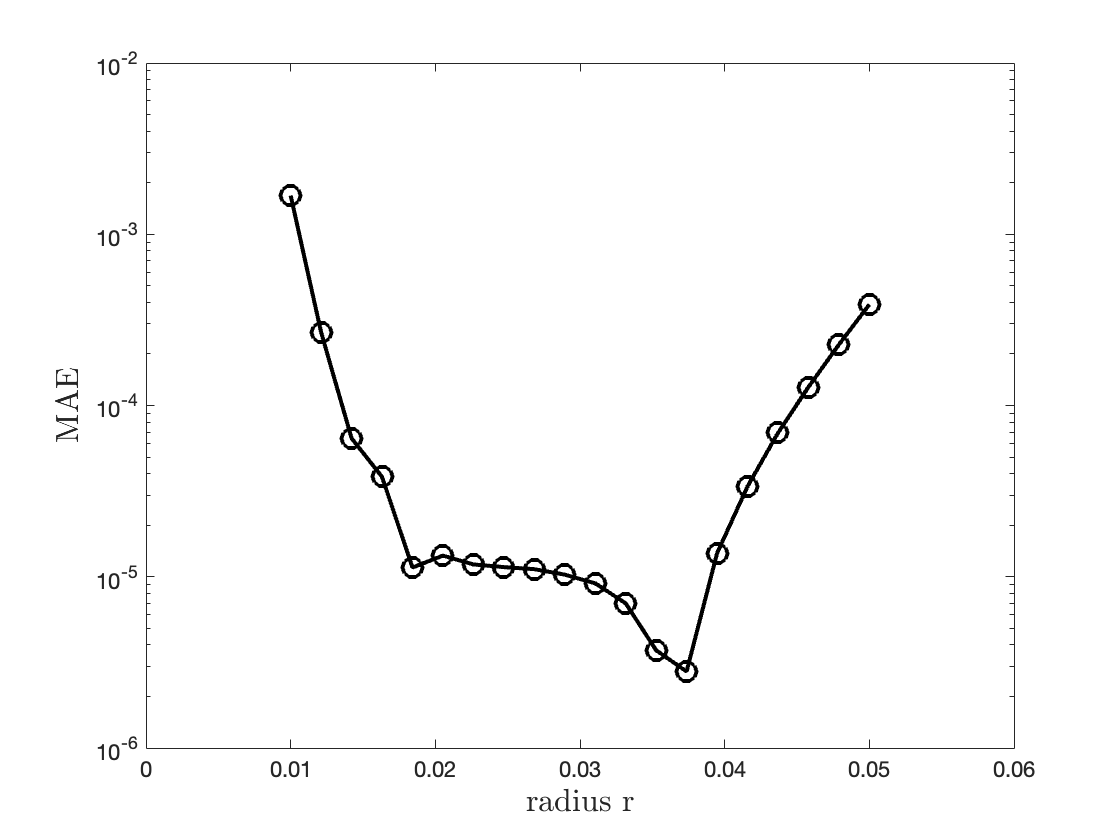}}
\caption{Quadratic Model: (a) Testing the influence of truncation order $N$ (b) Testing the influence of radius $r$}
\label{Fig.1}
\end{figure}

Figure \ref{Fig.1} summarizes these results in plots \ref{Fig.1a} and \ref{Fig.1b}. Test \ref{Fig.1a} shows that the exponential convergence of our approach with respect to truncation order, which is similar to the conclusion in conventional Carleman Linearization Method. Test \ref{Fig.1b} illustrates that the the error shows a "bowl shape" with respect to the radius $r$ i.e. $r$ cannot be too large or too small. 

\subsubsection{Cosine Square Model}\label{sec4.1.2}
\begin{figure}[H]
\centering  
\subfigure[Test of Truncation Order]{
\label{Fig.2a}
\includegraphics[width=0.45\textwidth]{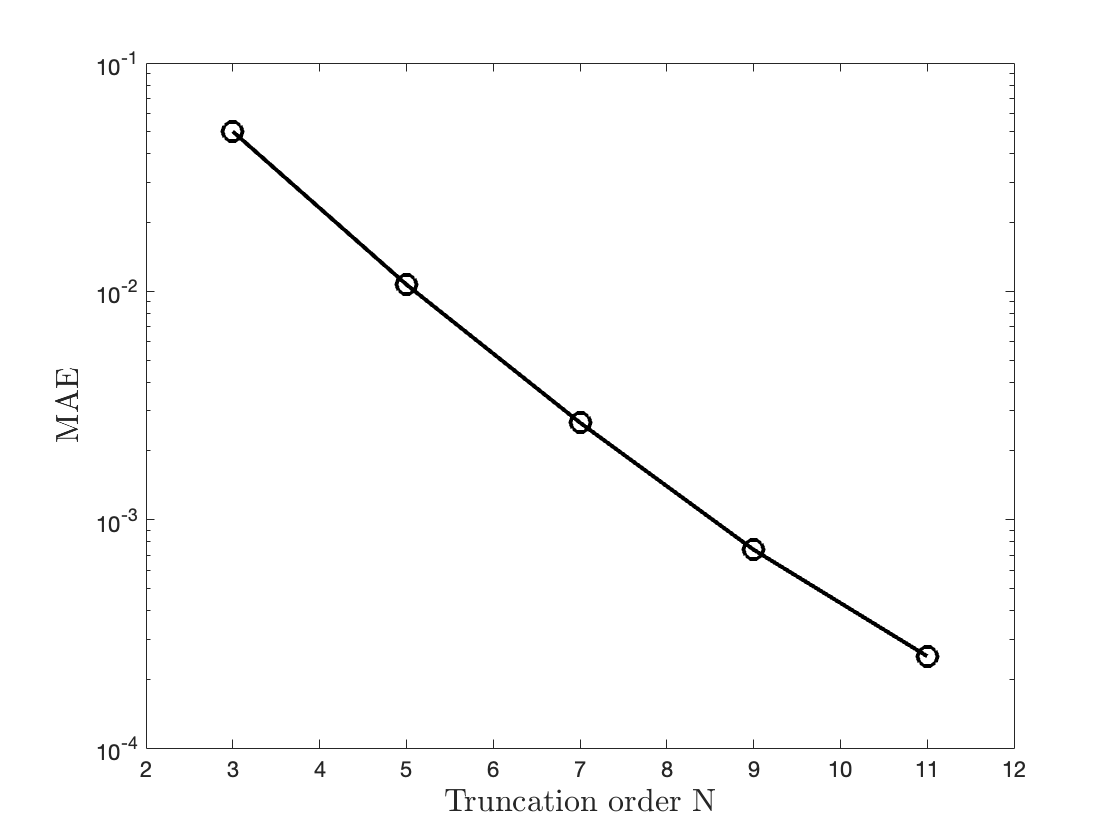}}
\subfigure[Test of Radius]{
\label{Fig.2b}
\includegraphics[width=0.45\textwidth]{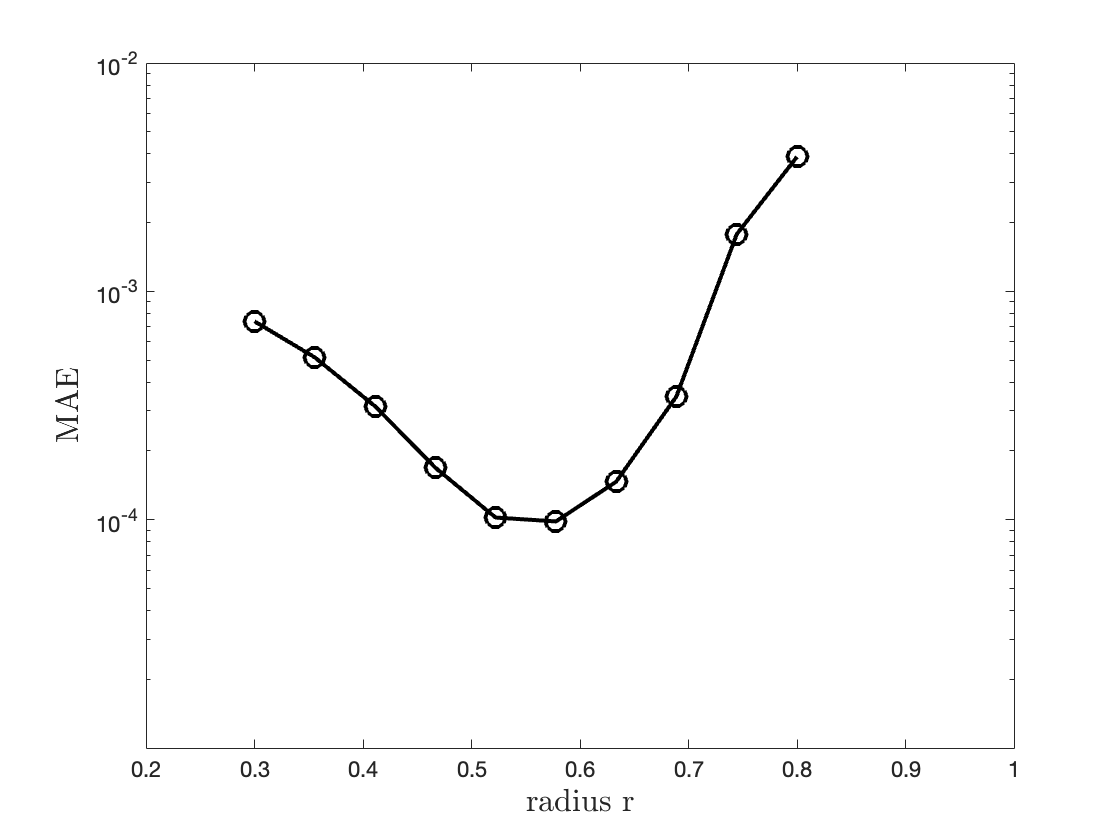}}
\caption{Cosine Square Model: (a) Testing Truncation Order $N$; (b) Testing the influence of radius $r$}
\label{Fig.2}
\end{figure}

This cosine square model is a synthetic model invented as a non-polynomial case for our demonstrative purposes. The governing equation is given as 
$$
\frac{dx}{dt}=\cos^2(x)
$$
And we set $x(0)=0.9$, $T=10$. Despite the nonlinearity, we still have a closed form solution $x(t)=\arctan(-0.5t+\tan(x_0))$. The three experiments are under the following parameters settings:
\begin{itemize}
  \item Test of Truncation Order N: $r=0.3$
  \item Test of radius r: $N=9$
\end{itemize}
Figure \ref{Fig.2} summarizes these results in plots \ref{Fig.2a} and \ref{Fig.2b}. Test \ref{Fig.2a} shows that the exponential convergence of our approach with
respect to truncation order, which is still similar to the conclusion in conventional Carleman Linearization
Method. Further, test \ref{Fig.2b} illustrates that the the error shows a ”V shape” with respect to the radius $r$ (i.e.
$r$ cannot be too large or too small as in the first example).

\subsubsection{Simple Pendulum Model}\label{sec4.1.3}
The simple pendulum is well studied model in science and engineering. The movement of the pendulum is described by a second order ordinary different equation,

$$
\frac{d^2\theta}{dt}=-\frac{g}{L}\sin(\theta)
$$

Where L is the length of the pendulum, $\theta$ is the displacement angle, and the parameter g is the gravity acceleration. This second order equation can be converted to a two dimensional first order ode system. As a notation we define $x_1=\theta$ and $x_2=\frac{d\theta}{dt}$. Also, we set $L=g=9.8$ for simplicity. Correspondingly,

\begin{align}
\frac{dx_1}{dt}&=x_2 \notag \\
\frac{dx_2}{dt}&=-\sin(x_1)    \notag
\end{align}

And we set $x(0)=(0.1,0.1)^T$, $T=10$. The three experiments are under the following parameters settings:
\begin{itemize}
  \item Test of solution: $r=(1,1)^T$;
  \item Test of Truncation Order $N$: $r=(1,1)^T$;
  \item Test of radius $r$: N=9
\end{itemize}

\begin{figure}[H]
\centering  
\subfigure[Test of Solution]{
\label{Fig.3a}
\includegraphics[width=0.45\textwidth]{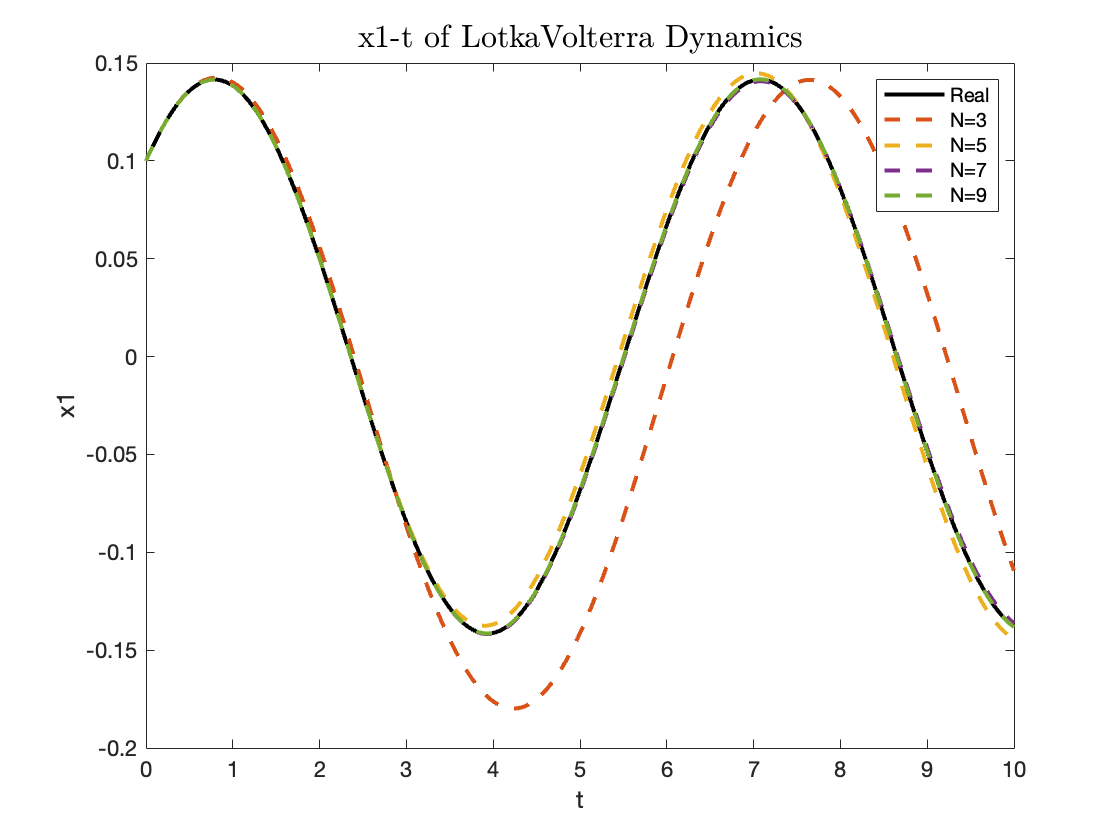}}
\subfigure[Test of Truncation Order]{\includegraphics[width=0.45\textwidth]{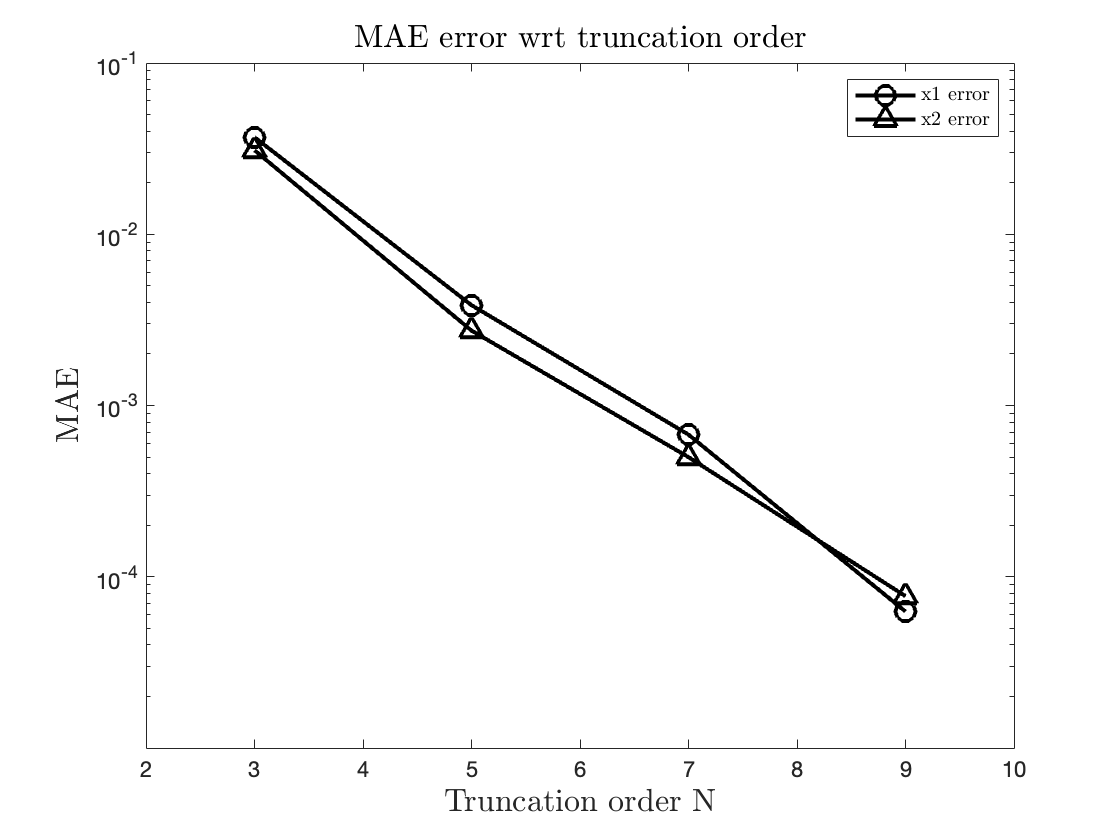}
\label{Fig.3b}}
\subfigure[Test of Radius]{
\label{Fig.3c}
\includegraphics[width=0.45\textwidth]{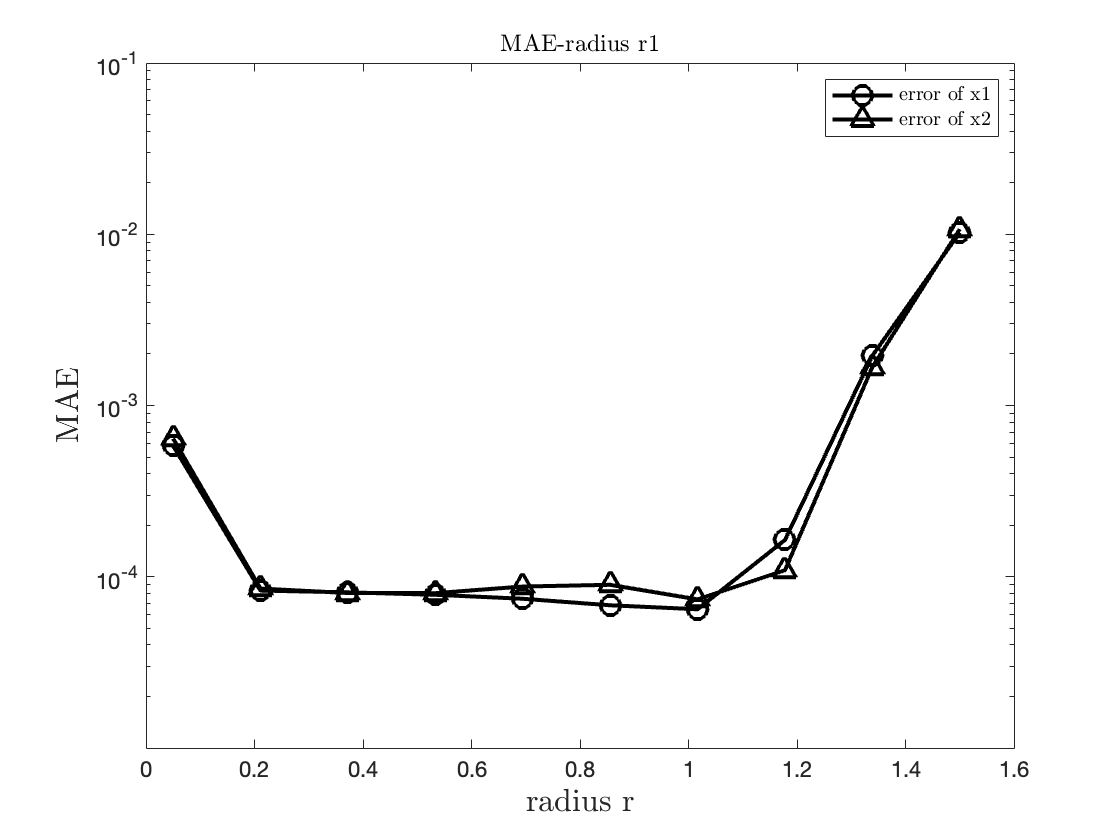}
\includegraphics[width=0.45\textwidth]{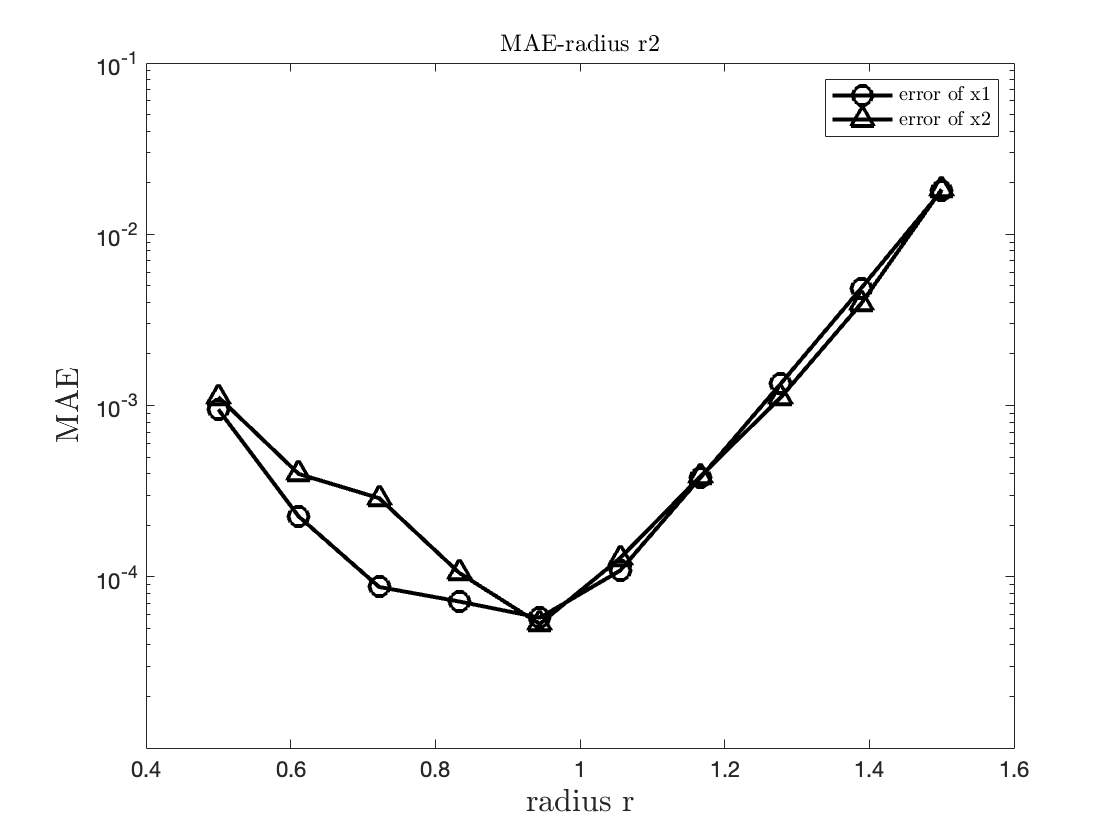}}
\caption{Simple Pendulum Model: (a) Testing Solution; (b) Testing Truncation Order $N$; (c) Testing the influence of radius $r$}
\label{Fig.3}
\end{figure}
Figure \ref{Fig.3} presents the numerical result of Simple Pendulum. Similar to previous result, the exponential convergence of error with respect to truncation order $N$ is observed from figure \ref{Fig.3b}. \ref{Fig.3c} indicate that $r_2=1$ might be the optimal setting and $r_1$ various from 0.2 to 1.2 give similar error result.
\subsubsection{Lotka-Volterra Model}\label{sec4.1.4}
The Lotka-Volterra model, also known as the predator-prey model\cite{ref29}, is a mathematical model used to describe the interactions between predators and prey in an ecosystem.The model describes the interaction between two fundamental populations: the prey population and the predator population. Here we defined as below

\begin{align}
\frac{dx_1}{dt}&=1.1 x_1-0.4x_1x_2 \notag \\
\frac{dx_2}{dt}&=0.1 x_1x_2-0.4 x_2    \notag
\end{align}

We set $x(0)=(5,5)^T$ and $T=10$. The parameters used in three different tests are as follows:
\begin{itemize}
  \item Test of solution: $r=(3,3)$;
  \item Test of Truncation Order $N$: $r=(3,3)^T$;
  \item Test of radius $r$: $N=13$
\end{itemize}

\begin{figure}[H]
\centering  
\subfigure[Test of Solution]{
\label{Fig.4a}
\includegraphics[width=0.45\textwidth]{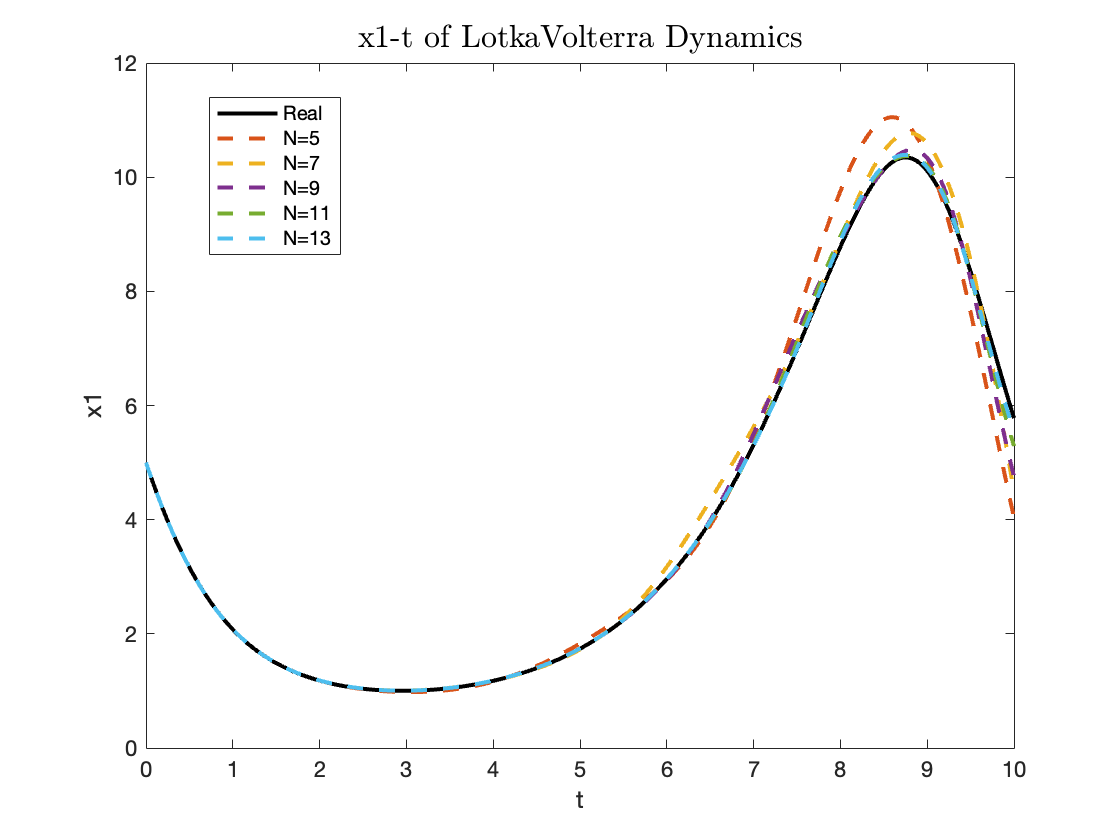}}
\subfigure[Test of Truncation Order]{
\label{Fig.4b}
\includegraphics[width=0.45\textwidth]{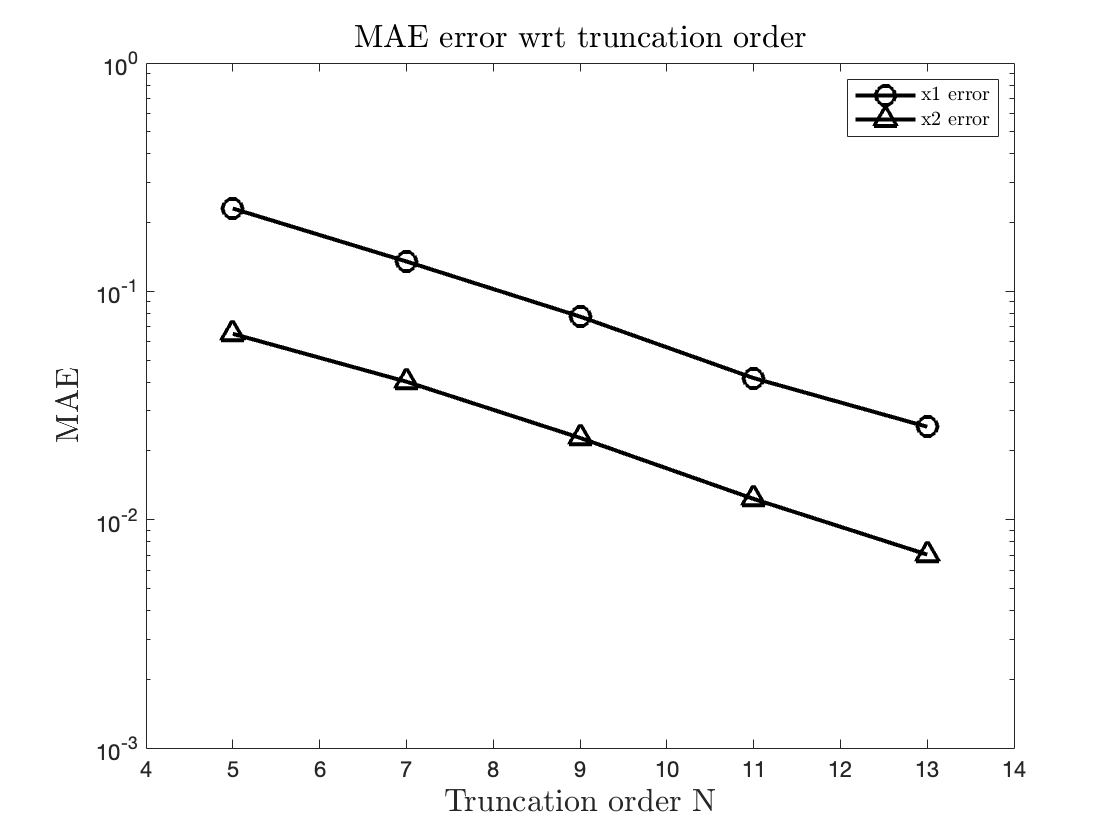}}
\subfigure[Test of Radius]{
\label{Fig.4c}
\includegraphics[width=0.45\textwidth]{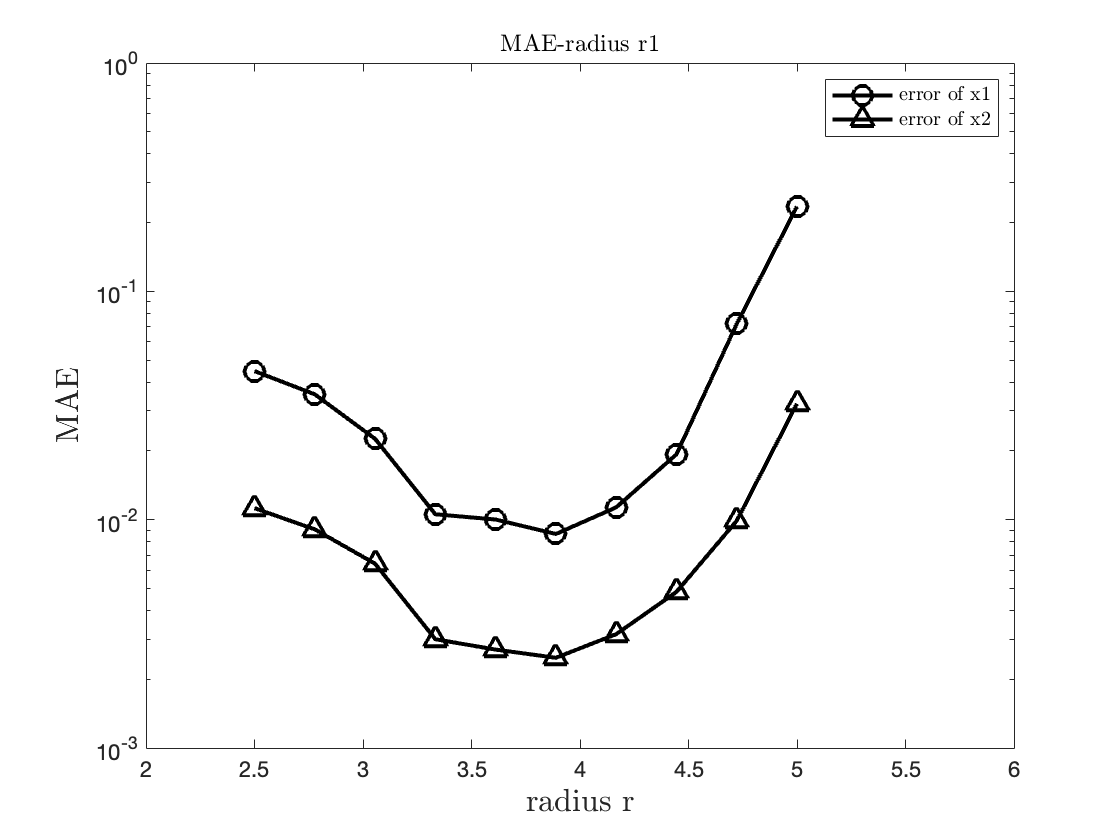}
\includegraphics[width=0.45\textwidth]{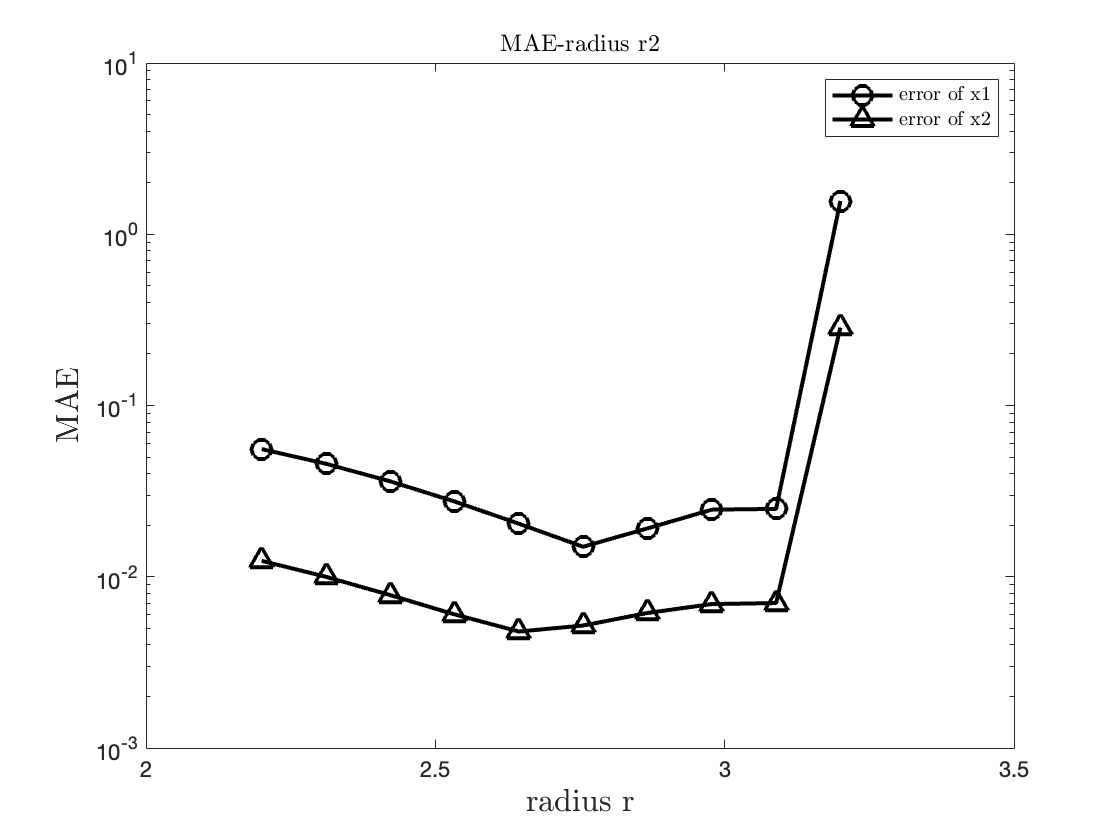}}
\caption{Lotka-Volterra Model: (a) Testing Solution; (b) Testing Truncation Order $N$; (c) Testing the influence of radius $r$}
\label{Fig.4}
\end{figure}
Figure \ref{Fig.4} presents the results of these tests for Lotka-Volterra Model. Again error decrease exponentially with respect to truncation order and "V shape" radius test results was observed. In particular, different from previous example, the results obtained between two different dimensions exhibit a consistent pattern, but the precision of the numerical outcomes for the second dimension is consistently higher in any given test.

\subsubsection{Kraichnan-Orszag Model}\label{sec4.1.5}
Kraichnan-Orszag Model was raised in \cite{ref26} for modeling fluid dynamics. This three dimensional nonlinear model is given by 
\begin{align}
\frac{dx_1}{dt}&=x_2x_3 \notag \\
\frac{dx_2}{dt}&=x_1x_3 \notag \\
\frac{dx_3}{dt}&=-2x_1x_2\notag
\end{align}
We set $x(0)=(0.1,-0.2,0.3)^T$ and $T=5$. In this experiment we apply the following parameter settings:
\begin{itemize}
  \item Test of solution: $r=(0.1,0.1,0.1)$,
  \item Test of Truncation Order $N$: $r=(0.1,0.1,0.1)^T$;
  \item Test of radius $r$: $N=9$
\end{itemize}

\begin{figure}[H]
\centering  
\subfigure[Test of Truncation Order]{
\label{Fig.5a}
\includegraphics[width=0.45\textwidth]{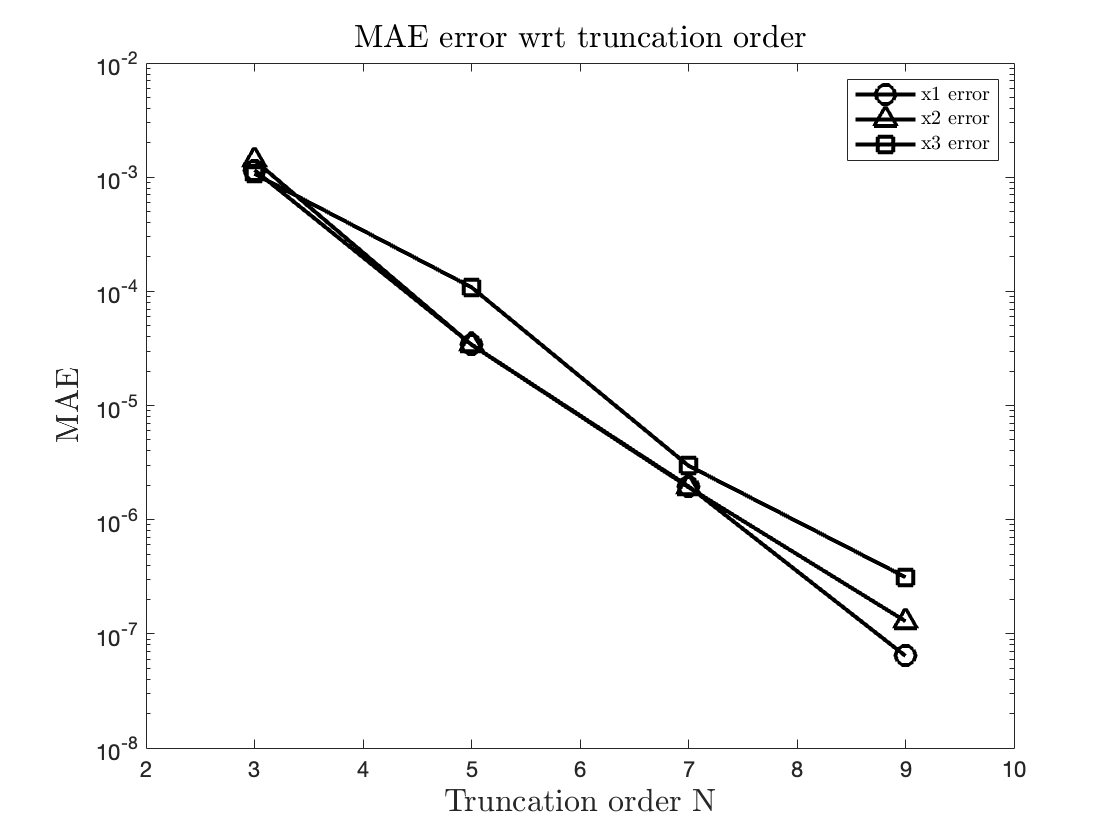}}
\subfigure[Test of Radius $r_1$]{
\label{Fig.5b}
\includegraphics[width=0.45\textwidth]{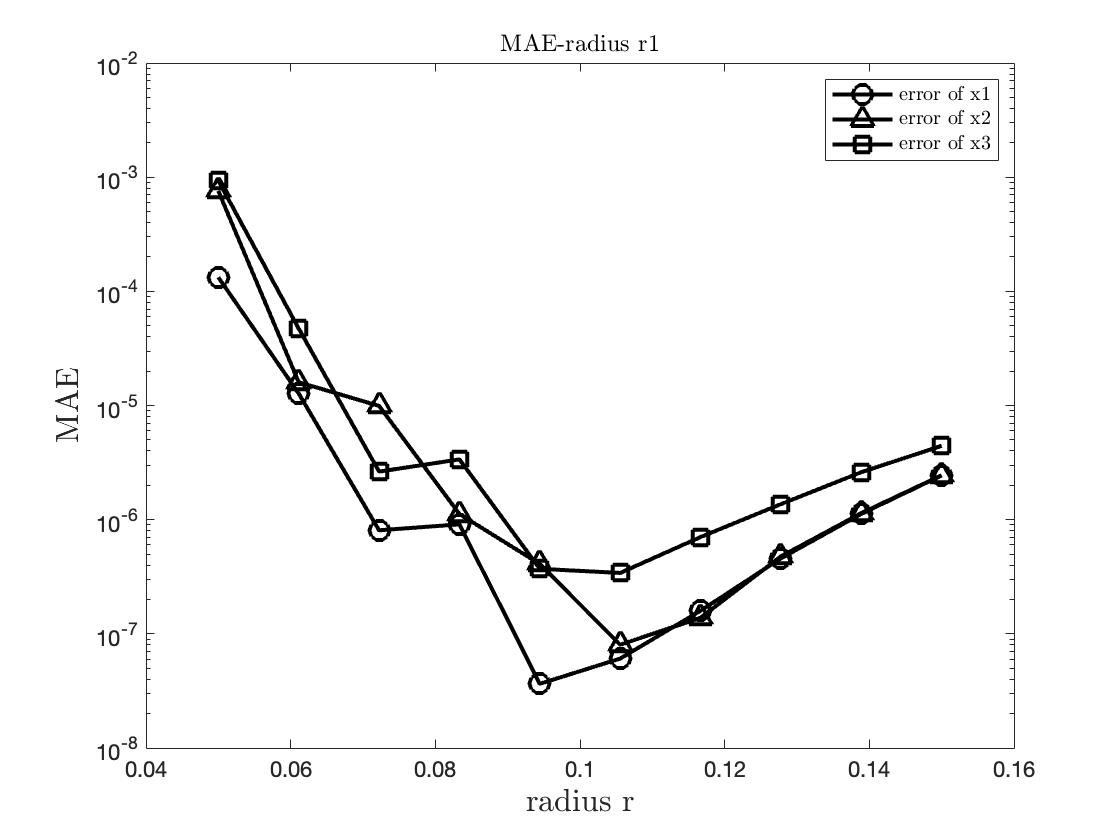}
}
\subfigure[Test of Radius $r_2$]{
\label{Fig.5c}\includegraphics[width=0.45\textwidth]{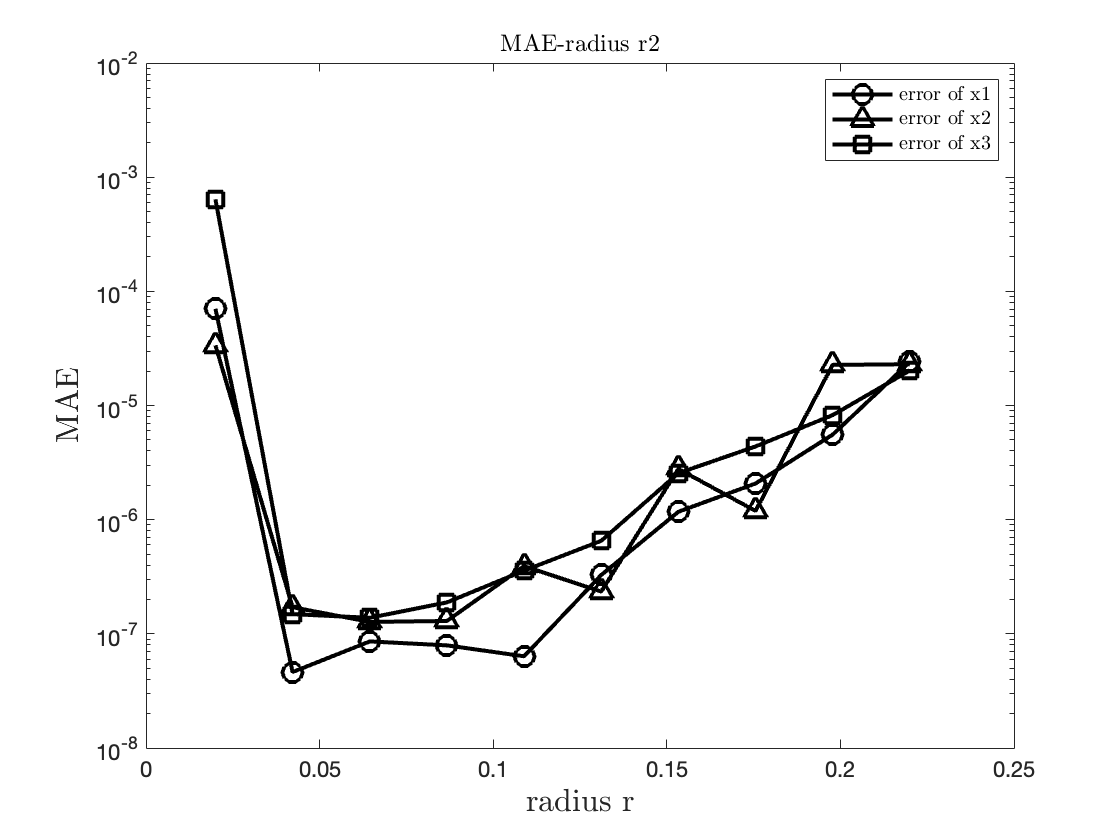}
}
\subfigure[Test of Radius $r_3$]{
\label{Fig.5d}\includegraphics[width=0.45\textwidth]{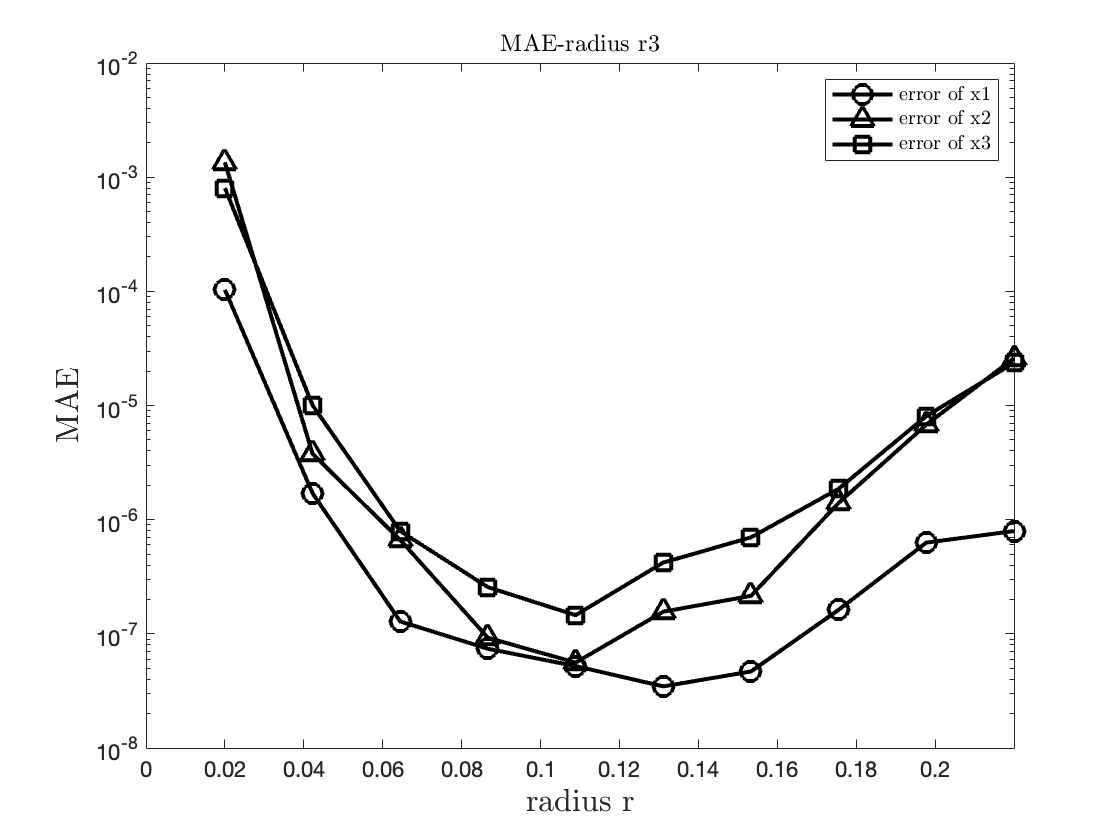}}

\caption{Kraichnan-Orszag Model: (a) Testing Truncation Order N; (b) Testing the influence of radius $r_1/r_2/r_3$}
\label{Fig.5}
\end{figure}
Figure \ref{Fig.5} presents the results of Kraichanan-Orszag Model. Different from previous example, a shorter time period is setting due to the strong oscillations of the Kraichnan-Orszag model. And based on observation of figure \ref{Fig.5b},\ref{Fig.5c}, and \ref{Fig.5d}, we can infer there might be an "best" $r$ choose around $r=(0.1,0.05,0.1)^T$
\subsection{Compared with Carleman Linearization}\label{sec4.2}
Following the construction procedure introduced in \autoref{sec3}, we apply the tensor product rule to construct collocation points. Consequently, the approximation matrix of the Koopman generator has a size of $N^d \times N^d$. However, for multi-dimensional dynamics, we can reuse the lifted matrix with different initial values. In the case of the Carleman linearization procedure, the matrix size will be $\sum_{i=1}^N d^i \times \sum_{i=1}^N d^i $, and we can approximate the size of the Carleman lifted matrix as $O(d^N)$. Therefore, when considering the matrix size of the derived linear system with a relatively small truncation order $N$, the matrix in our approach can be significantly smaller compared to the matrix in Carleman linearization. To further illustrate this point, \ref{tab.1} presents a matrix size comparison between our approach and Carleman linearization.
\setcounter{figure}{0}
\begin{figure}[H]
\centering  
\includegraphics[width=0.8\textwidth]{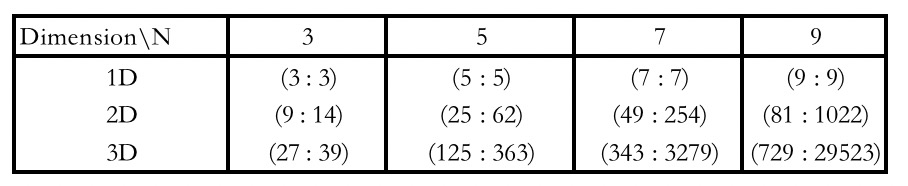}
\captionsetup{name={Table}}
\caption{Matrix Size Comparision (Koopman : Carleman)} 
\label{tab.1}
\end{figure}

\par While our approach generates $d$ linear systems with the same matrix but different initial values, which may not be as convenient as Carleman linearization, where only one linear system is produced, the significant difference in matrix sizes results in a notable increase in computational cost. As an illustrative example, we present the running time and error in \ref{tab.2} for a truncation order of $N$=9. It is evident that the computational time required for Carleman Linearization is slower than our approach, especially when $d$ = 2 or 3. In particular, for the Kraichnan-Orszag Model ($d$=3), the Carleman linearization procedure took 32.0546 seconds, which is considerably slower than our approach.
\par  Figure \ref{Fig.6} illustrates the error as a function of the truncation order $N$ and provides a comparison between Carleman linearization and our method. The error of Carleman linearization exhibits exponential decay, a phenomenon that has been substantiated by recent theoretical research \cite{ref27,ref28}. Our method leverages the Chebyshev node interpolation approach, which also enjoys exponential convergence guarantees and demonstrates similar exponential convergence in numerical experiments. For all the polynomial examples in \ref{Fig.6a}, \ref{Fig.6b}, and \ref{Fig.6c}, it is evident that our approach exhibits higher accuracy at the same truncation order, and the error decreases more rapidly as the truncation order is reduced.

Moreover, when applied to non-polynomial dynamical systems, our method showcases superior numerical accuracy and faster convergence rates. In Figure \ref{Fig.7}, we present a comparative analysis of the error-truncation order relationship in the context of non-polynomial dynamical systems. Both non-polynomial dynamics are subjected to 12th-order Taylor expansions to transform them into polynomial forms. It is evident that the Carleman method exhibits significantly lower accuracy in this scenario, with notably slower convergence. This limitation arises from the fact that, in non-polynomial cases, the relevant part of the Taylor expansion is confined to orders lower than the truncation order, resulting in substantial errors when employing relatively small expansion orders, as demonstrated in our study. Consequently, our method demonstrates heightened precision in non-polynomial cases. Even in the case of the 'simple pendulum' where $N$=3, our method initially exhibits lower accuracy. However, as $N$ increases to 9, the errors obtained by our method rapidly become dramatically smaller than those produced by the Carleman method.

\setcounter{figure}{1}
\begin{figure}[H]
\centering  
\includegraphics[width=0.8\textwidth]{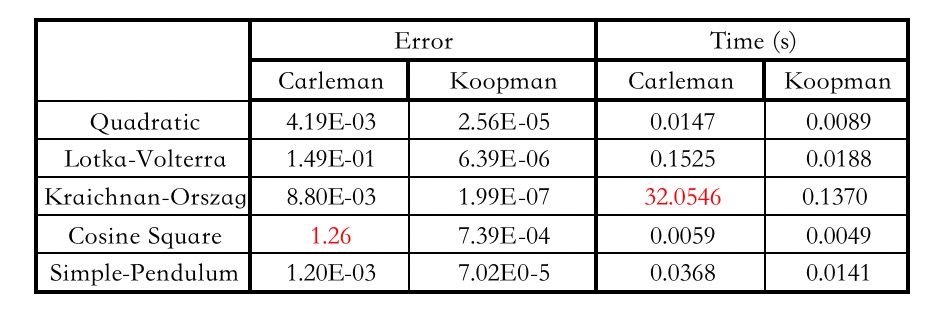}
\captionsetup{name={Table}}
\caption{Error and Time cost compared with Carleman Linearization ($N=9$)}
\label{tab.2}
\end{figure}

\setcounter{figure}{5}
\begin{figure}[H]
\centering  
\subfigure[Quadratic Model]{
\label{Fig.6a}
\includegraphics[width=0.3\textwidth]{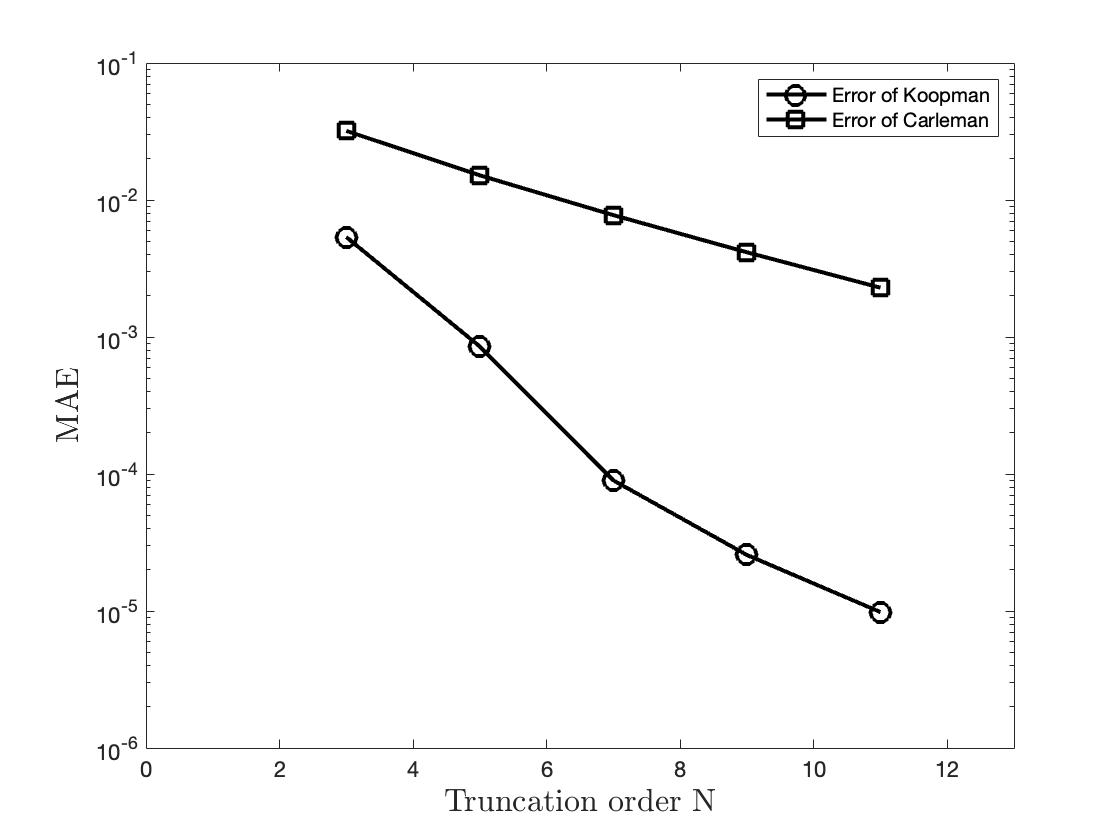}}
\subfigure[Lotka-Volterra Model]{
\label{Fig.6b}
\includegraphics[width=0.3\textwidth]{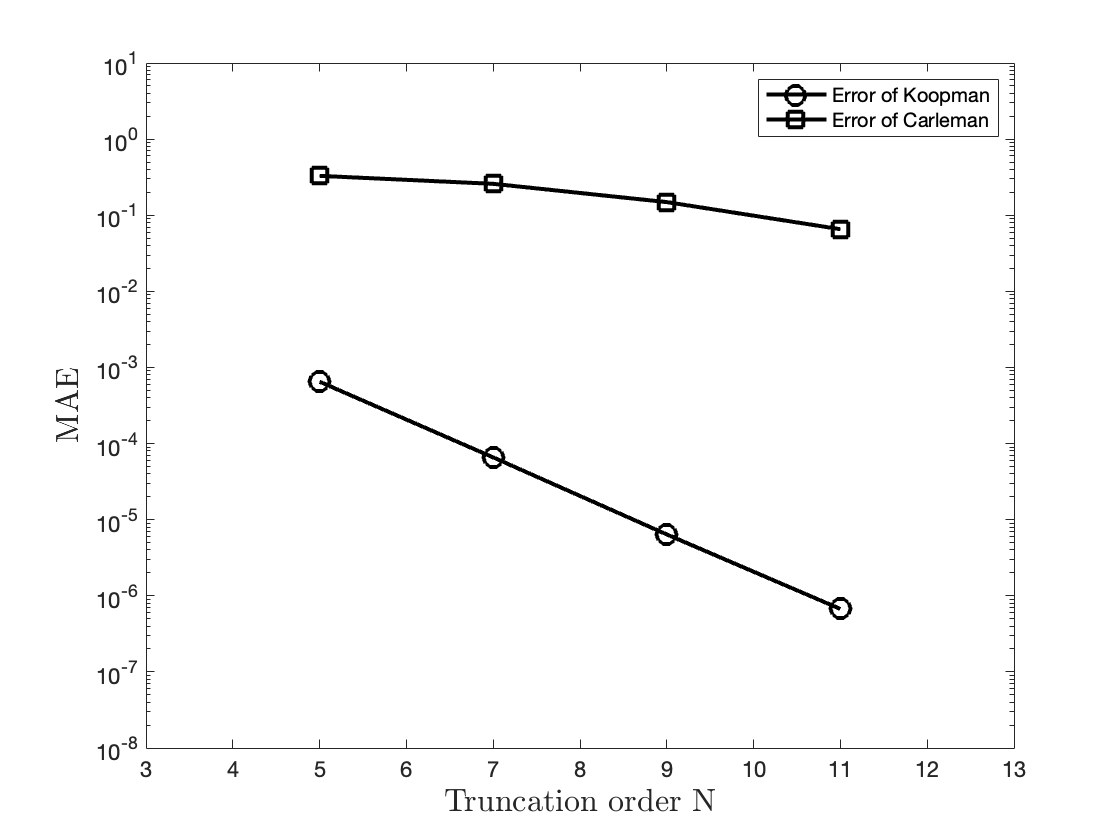}}
\subfigure[Kraichnan-Orszag Model]{
\label{Fig.6c}
\includegraphics[width=0.3\textwidth]{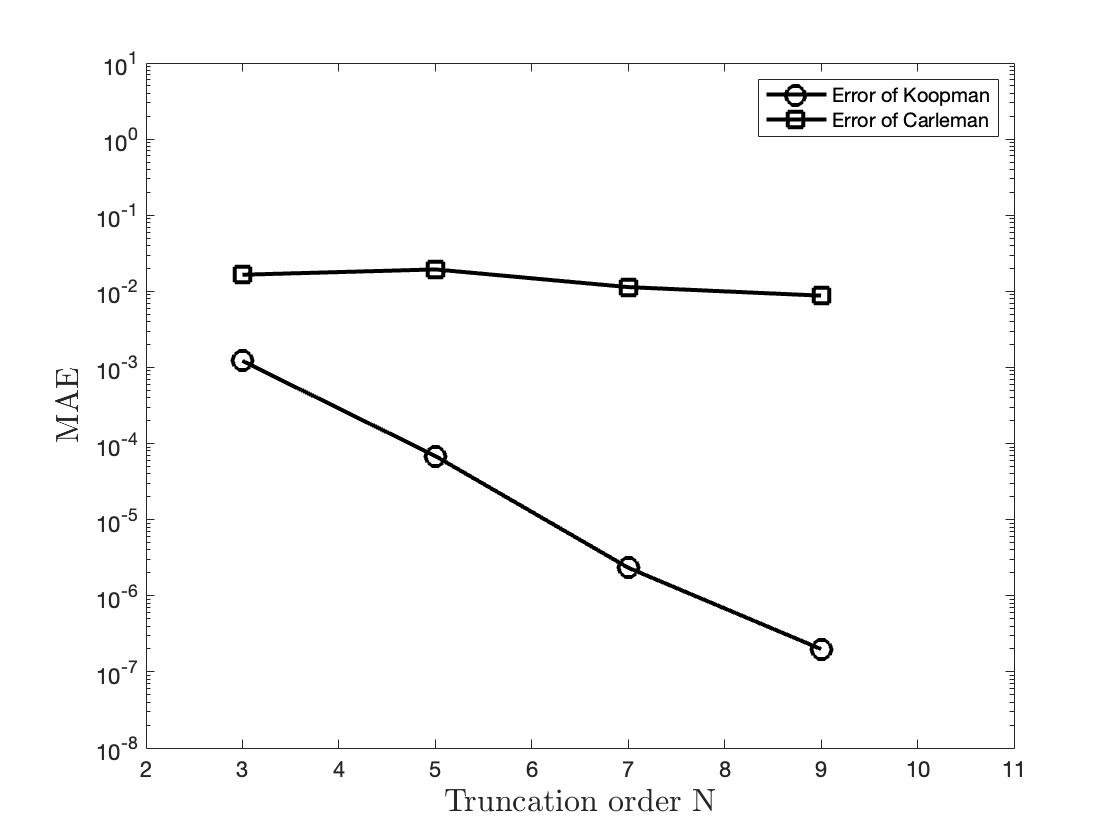}}
\caption{Error-Truncation Order compared with Carleman Linearization polynomial case (a) Quadratic Model; (b) Lotka-Volterra Model; (c) Kraichnan-Orszag Model}
\label{Fig.6}
\end{figure}
\begin{figure}[H]
\centering  

\subfigure[Cosine Square Model]{
\label{Fig.7a}
\includegraphics[width=0.45\textwidth]{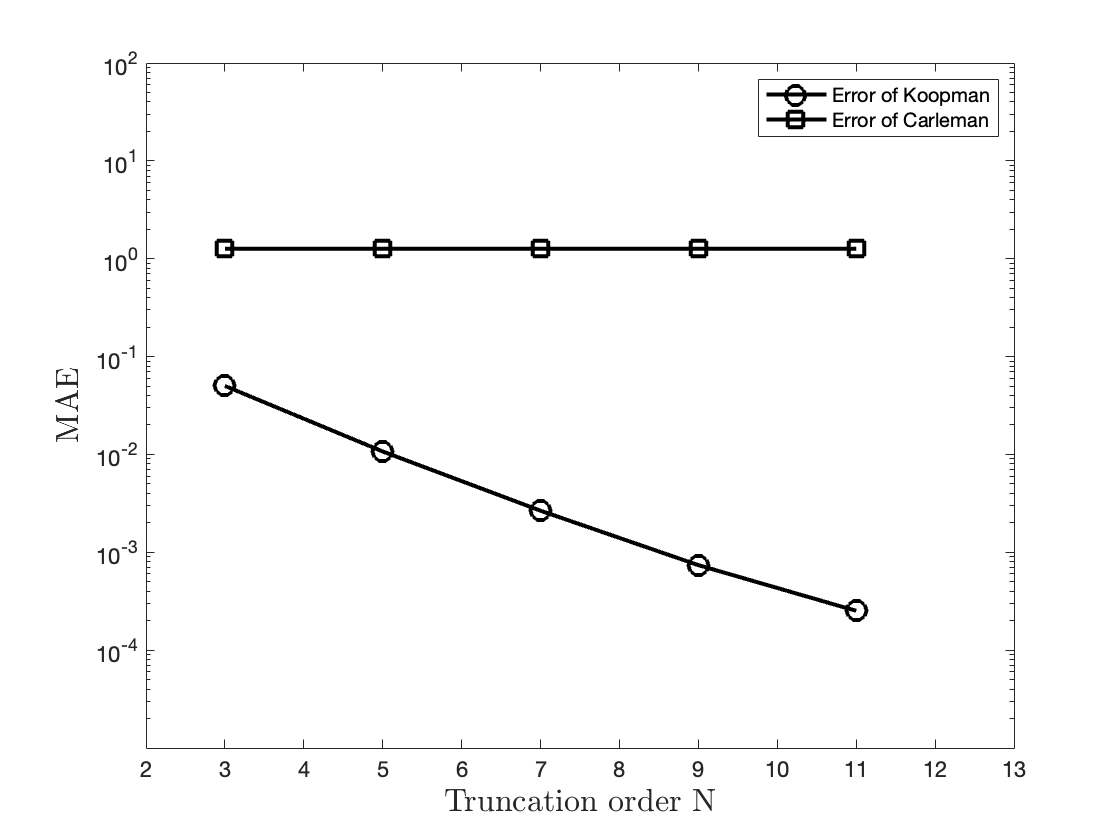}}
\subfigure[Simple Pendulum Model]{
\label{Fig.7b}
\includegraphics[width=0.45\textwidth]{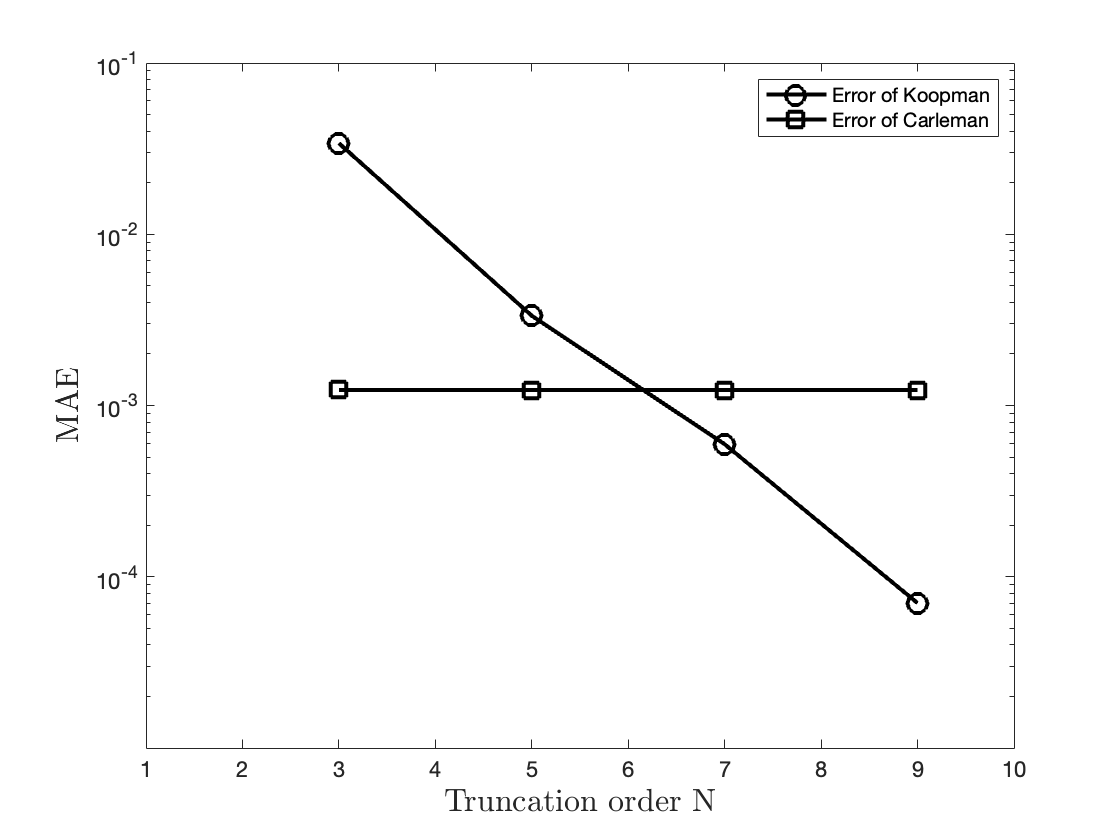}}
\caption{Error-Truncation Order compared with Carleman Linearization non-polynomial case (a) Cosine Square Model; (b) Simple Pendulum Model}
\label{Fig.7}
\end{figure}
With regard to limitations, our method's capacity to yield numerical results with smaller matrices, reduced computational time, and increased accuracy, all under the same truncation order, hinges upon a critical prerequisite—namely, the selection of an appropriate parameter $r$, a requirement not present in the Carleman Linearization method. As demonstrated in the results of the parameter $r$ testing discussed in \autoref{sec4.1}, deviating $r$ from the optimal setting point or range, either by making it excessively large or excessively small, results in gradually increasing errors.

Furthermore, when only the time-domain evolution differential equations of the dynamical system are available in advance, without access to spatial evolution information of the state, our method cannot capture the state's spatial evolution. Therefore, unless supplementary information is accessible regarding the spatial evolution of the dynamical system, allowing for a reasonable estimation of $r$, surpassing the precision of Carleman Linearization remains unattainable.

\section{Conclusion and Discussion}\label{sec5}
The Koopman Spectral Linearization method uses the Chebyshev Differentiation Matrix to explicitly derive the lifted matrix of nonlinear autonomous dynamical systems. It provides a finite-dimensional representation for the generator of the Koopman Operator with a polynomial basis. Therefore, similar to Carleman Linearization, our approach provides a linearization method for nonlinear dynamical systems. In each numerical experiment presented, our approach exhibits exponential convergence with respect to truncation order $N$, just like Carleman Linearization. Under the same truncation order, Koopman Spectral Linearization tends to exhibit significantly higher accuracy associated with lower time cost compared to Carleman Linearization, especially when the dynamics are not polynomial. Therefore, it is more suitable as an alternative linearization method where high-accuracy approximations are needed and $f$ is non-polynomial. Different from Carleman Linearization, which describes the evolution of the state itself, our approach is used to describe a scalar smooth observable, which is more flexible.
\par
However, despite its advantages, the Koopman Spectral Linearization method is not without its limitations. It's important to consider certain drawbacks when evaluating its applicability since a reasonable estimation of radius $r$ is necessary to obtain an accurate approximation of the eigenfunctions. Therefore if no additional information like range of states are given, its difficult to provide an accurate approximation. To further understand how this parameter $r$ influence the accuracy, a comprehensive numerical analysis will be included in the future work. Besides, the global interpolation method might overcome the limitation of estimate $r$.
\par
Regarding the matrix size as indicated in \autoref{sec4.2}, Koopman Spectral Method obtained a much smaller matrix with a even higher accuracy compared with Carleman Linearization. However the matrix size is exponential increase with respect to dimension d (i.e. $N^d$) since tensor product rule is applied. Therefore both Carleman Linearization (i.e. $O(d^N)$) and our approach might be not efficiency for relatively high-dimensional problem. A possible way to overcome this limitation is to adapt sparse grid methods to construct collocation points which has been founded success in \cite{ref17}. As pointed in Koopman Operator theory, either Carleman Linearization and our approach relies on tensor product rule, therefore cannot deal with high dimension problems. By combining with sparse grid technique, our approach is possible to address much higher-dimensional problems\cite{ref17}.
\par
Furthermore, an interesting potential application about our approach is to address efficient quantum algorithm for nonlinear ODE by combination of quantum linear system algorithm to solve the linearized system \cite{ref7}. Especially our approach only require the middle element of the output vector which might significantly improve the process of extract quantum information.

\newpage
\printbibliography
\end{document}